\documentclass[11pt]{amsart}
\usepackage{palatino}
\usepackage{amsfonts}
\usepackage{amssymb}
\usepackage{amscd}

\theoremstyle{plain}

\newtheorem{theorem}{Theorem}[section]
\newtheorem{proposition}[theorem]{Proposition}

\newtheorem{lemma}[theorem]{Lemma}
\theoremstyle{definition}
\newtheorem{definition}[theorem]{Definition}
\newtheorem{remark}[theorem]{Remark}

\newtheorem{conjecture}[theorem]{Conjecture}
\newtheorem{conjecture/question}[theorem]{Conjecture/Question}

\newtheorem{remark/definition}[theorem]{Remark/Definition}
\newtheorem{terminology/notation}[theorem]{Terminology/Notation}

\newcommand{\marginlabel}[1]%
  {\mbox{}\marginpar{\raggedleft\hspace{0pt}\bfseries\sf#1}}

\def\rmapdown#1{\Big\downarrow
   \rlap{$\vcenter{\hbox{$\scriptstyle#1$}}$ }}

\def\QQ{{\mathbb Q}}
\def\PP{{\textbf P}}
\def\OO{\mathcal{O}}
\def\cN{\mathcal{N}}

\def\cB{\mathcal{B}}
\def\cA{\mathcal{A}}
\def\F{\mathcal{F}}
\def\E{\mathcal{E}}
\def\G{\mathcal{G}}
\def\K{\mathcal{K}}
\def\L{\mathcal{L}}
\def\I{\mathcal{I}}

\def\cM{\mathcal{M}}

\def\cU{\mathcal{U}}
\def\cC{\mathcal{C}}

\def\H{\mathcal{H}}
\def\Pic0{{\rm Pic}^0(X)}

\def\mm{\overline{\mathcal{M}}}
\def\kk{\overline{\mathcal{K}}}
\def\zz{\overline{\mathcal{Z}}}

\theoremstyle{remark}

\begin{document}

\title{\bf Syzygies of curves and the effective cone of $\mm_g$}

\author[G. Farkas]{Gavril Farkas}
\address{Department of Mathematics, University of Texas,
Austin, TX 78712} \email{{\tt gfarkas@math.utexas.edu}}
\thanks{Research  partially supported by the NSF Grant DMS-0450670 and by the Sloan Foundation}

\maketitle


\section{Introduction}

The aim of this paper is to describe a systematic way of
constructing effective divisors on $\mm_g$ having exceptionally
small slope. In particular, these divisors provide a string of
counterexamples to the Harris-Morrison Slope Conjecture (cf.
\cite{HMo}). In a previous paper \cite{FP}, we showed that the
divisor $\kk_{10}$ on $\mm_{10}$ consisting of sections of $K3$
surfaces contradicts the Slope Conjecture on $\mm_{10}$. Since the
moduli spaces $\mm_g$ are known to behave erratically for small
$g$ and since the condition that a curve of genus $g$ lie on a
$K3$ surface is divisorial only for $g=10$, the question remained
whether $\kk_{10}$ is an isolated example or the first in a series
of counterexamples. Here we prove that any effective divisor on
$\mm_g$ consisting of curves satisfying a Green-Lazarsfeld syzygy
type condition for a linear system residual to a pencil of minimal
degree, violates the Slope Conjecture. A consequence of the
existence of these effective divisors is that various moduli
spaces $\mm_{g, n}$ with $g\leq 22$, are proved to be of general
type.

We recall that the slope $s(D)$ of an  effective divisor $D$ on
$\mm_g$ is defined as the smallest rational number $a/b \geq 0$
such that the divisor class $a
\lambda-b(\delta_0+\cdots+\delta_{[g/2]})-[D]$ is an effective
combination of boundary divisors. The Slope Conjecture predicts
that $s(D)\geq 6+12/(g+1)$ for all effective divisors $D$ on
$\mm_g$ (cf. \cite{HMo}). More generally, the question of finding
a good lower bound for the \emph{slope of $\mm_g$}
 $$s_g:=\mbox{inf}\{s(D): D\in
\mbox{Eff}(\mm_g)\}$$ is of great interest for a variety of
reasons, for instance it would provide a new geometric solution to
the Schottky problem. In a different direction, since
$s(K_{\mm_g})=13/2$ (cf. \cite{HM}), to prove that $\mm_g$ is of
general type it suffices to exhibit a single effective divisor $D$
on $\mm_g$ of slope $s(D)<13/2$.

The Slope Conjecture is true for all $\mm_g$ with $g\leq 9$ but in
\cite{FP} we proved that on $\mm_{10}$, we have the equality
$s(\kk_{10})=7<6+12/11$ (in fact $\kk_{10}$ is the only effective
divisor on $\mm_{10}$ having slope $<6+12/11$). In \cite{FP} we
also showed that $\K_{10}$ has four incarnations as a geometric
subvariety of $\cM_{10}$. In particular, $\K_{10}$ can be thought
of as either

\noindent
(1) the locus of curves $[C]\in \cM_{10}$ for which
the rank $2$ Mukai type Brill-Noether locus $$\mbox{SU}_2(C, K_C,
6):=\{E\in SU_2(C, K_C): h^0(C, E)\geq 7\}$$ is not equal to $\emptyset$,
or

\noindent
(2) the locus of curves $[C]\in \cM_{10}$ carrying a pencil
$A\in W^1_6(C)$ such that the multiplication map
$\mu_A:\mbox{Sym}^2 H^0(C, K_C\otimes A^{\vee})\rightarrow H^0(C,
(K_C\otimes A^{\vee})^{\otimes 2})$ is not surjective.

The geometric conditions (1) and (2), unlike the original
definition of $\K_{10}$, can be extended to other genera. We fix
$g=2k-2$ and denote by $\sigma:\mathfrak G^1_k\rightarrow \cM_g$
the  Hurwitz scheme of $k$-sheeted coverings of $\PP^1$ of genus $g$
parametrizing pairs $(C, A)$ with $A\in W^1_k(C)$.
For each $i\geq 0$ we introduce the cycle $\mathcal{U}_{g, i}$
consisting of pairs $(C, A)\in \mathfrak G^1_k$ such that $K_C\otimes A^{\vee}$ fails the
Green-Lazarsfeld property $(N_i)$. By setting $\mathcal{Z}_{g, i}:=\sigma (\mathcal{U}_{g, i})$ we obtain
an induced geometric stratification of $\cM_g$
$$\mathcal{Z}_{g, 0}\subset \mathcal{Z}_{g, 1}\subset \ldots \subset
\mathcal{Z}_{g, i}\subset \ldots \subset \mathcal{M}_{g}.$$
If for  each $(C, A)\in \mathfrak G^1_k$ we consider the
 map $C\stackrel{|K_C\otimes A^{\vee}|}\longrightarrow \PP^{k-2}$ induced by
 the residual linear system, we can define two vector bundles $\cA$ and $\cB$ on $\mathfrak G^1_k$
 such that$$\cA(C, A)=H^0(\Omega_{\PP^{k-2}}^i(i+2)) \mbox{ and} \ \cB(C, A)=H^0(\Omega_{\PP^{k-2}}^i(i+2)\otimes \OO_C).$$
  For $g=6i+10$ (hence $k=3i+6$), it turns out that $\mbox{rank}(\cA)=\mbox{rank}(\cB)$ and $\mathcal{U}_{g, i}$ is the
  degeneracy locus of the natural vector bundle map $\phi: \cA\rightarrow \cB$, that is, $\mathcal{Z}_{g, i}$ is a virtual divisor on $\cM_g$, with virtual class $\sigma_* c_1(\cB-\cA)$.
The main result of this article is the computation of the compactification inside $\mm_g$ of this degeneracy locus (see Theorem \ref{formula} for a
 more precise statement):
\vskip 6pt \noindent \textbf{Theorem A. }\emph{ If
$\sigma:\overline{\mathfrak G}^1_{3i+6} \rightarrow \mm_{6i+10}$ is the
compactification of the Hurwitz stack by limit linear series,
then there is a  natural extension of the vector bundle map
$\phi:\cA\rightarrow \cB$ over $\overline{\mathfrak G}^1_{3i+6}$ such
that
 $\overline{\mathcal{Z}}_{g, i}$ is the degeneracy locus of $\phi$. Moreover the class of the
 pushforward to $\mm_g$ of the virtual degeneracy locus of $\phi$
 is given by
$$\sigma_*(c_1(\cB)-c_1(\cA))\equiv a\ \lambda-b_0\ \delta_0-b_1\
\delta_1-\cdots -b_{3i+5}\delta_{3i+5},$$
where  $a, b_0, \ldots, b_{3i+5}$ are explicitly determined coefficients such that  $b_j\geq b_0$ for $j\geq 1$ and
$$\frac{a}{b_0}=\frac{3(4i+7)(6i^2+19i+12)}{(i+2)(12i^2+31i+18)}<6+\frac{12}{g+1}.$$}

\noindent \textbf{Corollary.} \emph{If the degeneracy locus $\mathcal{Z}_{6i+10, i}$ is an actual divisor on $\cM_{6i+10}$, then we have that
$$s(\overline{\mathcal{Z}}_{6i+10, i})=\frac{3(4i+7)(6i^2+19i+12)}{(i+2)(12i^2+31i+18)}<6+\frac{12}{g+1},$$
thus contradicting the Slope Conjecture on $\mm_{6i+10}$.}

\vskip 5pt
The idea of the proof is to define a whole host of vector bundles $\G_{a, b}$ and $\H_{a, b}$ over $\mathfrak G^1_k$ for $a, b\geq 0$, having fibres
$$\G_{a, b}(C, A)=H^0(\Omega_{\PP^{k-2}}^a(a+b)\otimes \OO_C)\ \mbox{ and } \ \H_{a, b}(C, A)=H^0(\Omega_{\PP^{k-2}}^a (a+b)).
$$
 These bundles are related to one another by certain exact sequences (\ref{gi}) and (\ref{sym}) over $\mathfrak G^1_k$ (see
 Section 3). After an analysis over each boundary divisor $\sigma^{-1}(\Delta_j)$, where $0\leq j\leq 3i+5$,
 we find in Section 3 a unique way of extending $\G_{a, b}$ and $\H_{a, b}$ to
  vector bundles over $\overline{\mathfrak G}^1_k$ such
that (\ref{gi}) and (\ref{sym}) continue to make sense and be
exact and then use this to compute the Chern numbers of
$\cA=\H_{i, 2}$ and $\cB=\G_{i, 2}$.

We expect $\mathcal{Z}_{6i+10, i}$ to be always a divisor on
$\cM_{6i+10}$ but we are able to check this only for small $i$. To
verify this in general  one would have to prove that if $[C]\in
\cM_{6i+10}$ is a general curve, then one (or equivalently all) of
the finitely many  linear systems  $\mathfrak
g^{3i+4}_{9i+12}=K_C(-\mathfrak g^1_{3i+6})$ satisfies property
$(N_i)$. This is a direct generalization of Green's Conjecture to
the case of curves and line bundles $L$ with $h^1(L)=2$ (see Section
2 for more on this analogy).

For $g=10$ we recover of course the results from \cite{FP}. In the
 next two cases, $g=16$ and $22$  we have complete results:

\begin{theorem}\label{g16}
The following subvariety of $\cM_{16}$
$$\mathcal{Z}_{16, 1}:=\{[C]\in \cM_{16}: \exists L\in W^7_{21}(C) \mbox{ such that }C\stackrel{|L|}\hookrightarrow \PP^7 \mbox{ is not cut out by quadrics}\},$$
is an effective divisor, and the class of its compactification
is given by the formula
$$\zz_{16, 1}\equiv 286\bigl(407 \lambda-61\ \delta_0- 325\ \delta_1-b_2\ \delta_2-\cdots -b_8\ \delta_8\bigr),$$
where  $b_j\geq b_1$ for all $2\leq j\leq 8$. In particular
$s(\zz_{16, 1})=407/61 =6.6721...<6+12/17=6.705...$, hence $\zz_{16,
1}$ provides a counterexample to the Slope Conjecture on $\mm_{16}$.
\end{theorem}

In a similar manner we have the following example of a geometric divisor
on $\mm_{22}$ of very small slope:

\begin{theorem}\label{g22}
The following subvariety of $\cM_{22}$
$$\mathcal{Z}_{22, 2}:=\{[C]\in \cM_{22}: \exists L\in W^{10}_{30}(C) \mbox{ such that } C\stackrel{|L|}\hookrightarrow \PP^{10} \mbox{ fails property } (N_2)\},$$
is an effective divisor on $\cM_{22}$ and the class of its
compactification is given by the formula $$\zz_{22, 2}\equiv
25194\bigl(1665\ \lambda-256\  \delta_0- 1407\
\delta_1-b_2\delta_2-\cdots -b_{11}\  \delta_{11}\bigr),$$ where
$b_j\geq b_1$ for all $2\leq j\leq 11$. In particular $s(\zz_{22,
2})=1665/256=6.50390...<6+12/23=6. 52173...$, and $\zz_{22, 2}$
gives a counterexample to the Slope Conjecture on $\mm_{22}$.
\end{theorem}

We note that while Theorem A gives a sequence of (virtual) examples
of divisors contradicting the Slope Conjecture, for any given $g$
one can construct other (actual) divisors on $\mm_g$ of slope $<
6+12/(g+1)$ also defined  in terms of syzygies. For instance, the
locus consisting of curves $[C]\in \mm_{20}$ for which there exists
$L\in W^4_{20}(C)$ such that $C\stackrel{|L|}\hookrightarrow \PP^4$
is not cut out by \emph{quartics} (or equivalently, the map
$I_4(L)\otimes H^0(L)\rightarrow I_5(L)$ is not bijective), is a
divisor on $\mm_{20}$ violating the Slope Conjecture.

We record the following consequence of Theorems \ref{g16} and
\ref{g22}:

\begin{theorem}\label{mgn}
\noindent $(1)$ The moduli space of $n$-pointed curves $\mm_{22,
n}$ is of general type for all $n\geq 2$.

\noindent  $(2)$ The moduli space $\mm_{21, n}$ is of general type
for all $n\geq 4$.

\noindent  $(3)$ The moduli space $\mm_{16, n}$ is of general type
for all $n\geq 9$.

\noindent $(4)$ The moduli space $\mm_{20, n}$ is of general type for $n\geq 6$.
\end{theorem}

It is also possible to give examples of divisors of small slope on the
moduli spaces $\mm_{g, n}$ of $n$-pointed stable curves. We have
only pursued this for genera $g\leq 22$ with the goal of proving
that various moduli spaces $\mm_{g, n}$ are of general type (this
is automatic for $g\geq 23$). One of our examples is the
following:

\begin{theorem}\label{m14}
The following subvariety of $\mm_{14, 1}$
$$\mathcal{Z}_{14, 0}^1:=\{[C, p]\in \cM_{14, 1}:\exists L\in W^6_{18}(C) \mbox{ such that } C\stackrel{|L(-p)|}\longrightarrow \PP^5 \mbox{ fails property }(N_0)\}$$
is an effective divisor on $\cM_{14, 1}$. The class of its
compactification $\zz_{14, 0}^1$ is given by the formula
$$\zz_{14, 0}^1\equiv 33\bigl(237\ \lambda+14\ \psi-35\ \delta_0-169 \delta_1- 183 \delta_{13}- b_2\delta_2-\cdots
-b_{12}\delta_{12}),$$ where $b_j\geq 15+27j$ for $3\leq j\leq
12$, $j\neq 4$,  $b_2\geq 325$ and  $b_4\geq 271.$ In particular
$[\zz_{14, 0}^1]$ lies outside the cone of $\rm{Pic}$$(\mm_{14,
1})$ spanned by pullbacks of effective divisors from $\mm_{14}$,
the boundary divisors $\delta_0, \ldots, \delta_{13}$ and the
Weierstrass divisor $\overline{\mathcal{W}}$ on $\mm_{14, 1}$.
\end{theorem}

An equivalent formulation of the Slope Conjecture on $\mm_g$ for
$g$ such that $g+1$ is composite (which we learned from S. Keel),
is to say that the Brill-Noether divisors $\mm_{g, d}^r$
consisting of curves with a $\mathfrak g^r_d$ when
$g-(r+1)(g-d+r)=-1$, lie on a face of the effective cone
$\mbox{Eff}(\mm_g)$. It is well-known that $s(\mm_{g,
d}^r)=6+12/(g+1)$ (cf. \cite{EH3}), so the Slope Conjecture
singles out these divisors as being of minimal slope. One can ask
a similar question on $\mm_{g, 1}$. Are the pullbacks
$\pi^*(\mm_{g, d}^r)$ of the Brill-Noether divisors from $\mm_g$
and the Weierstrass divisor $\overline{\mathcal{W}}:=\{[C, p]\in
\mm_{g, 1}: p\in C \mbox{ is a Weierstrass point}\}$ on a face of
the effective cone $\mbox{Eff}(\mm_{g, 1})$? The question makes
sense especially since in \cite{EH2} it is proved that the class
of any generalized Brill-Noether divisor on $\mm_{g, 1}$ (that is,
any codimension $1$  locus of curves $(C, p)\in \mm_{g, 1}$ having
a linear series with special ramification at $p$) lies inside the
cone of $\mbox{Pic}(\mm_{g, 1})$ spanned by by $[\pi^*(\mm_{g,
d}^r)]$ and $[\overline{\mathcal{W}}]$. Theorem \ref{m14} shows
that at least for $g\geq 14$ the answer to the question raised
above is emphatically negative.

\textbf{Acknowledgments:} I  had useful conversations with many
people on subjects related to this paper. I especially benefitted
from discussions with
 J. Harris, S. Keel, D.
Khosla and M. Popa.
\section{Syzygies of algebraic curves}

In this paragraph we review a few facts about the resolution of
the ideal of a curve embedded in a projective space. As a general
reference for syzygies and Koszul cohomology
 we recommend \cite{L} or \cite{GL}.

 Suppose that $C$ is a smooth curve of genus $g$ and $L$ is
a very ample line bundle on $C$ giving an embedding $C\rightarrow
\PP^r=\PP(V)$, where $V=H^0(L)$. We denote by $I_{C/\PP^r}$ the
ideal of $C$ in $\PP^r$ and consider its minimal resolution of
free $S=\mbox{Sym}(V)$-modules
$$0\rightarrow F_{r+1}\rightarrow \cdots \rightarrow F_2\rightarrow F_1\rightarrow
I_{C/\PP^r}\rightarrow 0.$$ Then one can write $F_j=\oplus_{l\in
\mathbb Z} S(-j-l)^{b_{j,l}(C)}$, where
$b_{j,l}(C)=\mbox{dim}_{\mathbb C}\mbox{Tor}_j^S\bigl(S(C),
\mathbb C \bigr)_{j+l}$ is the graded Betti number of $C$ that
comes on the $l$-th row and $j$-th column in the Betti diagram of
$C$.  Following Green and Lazarsfeld we say that the pair $(C, L)$
satisfy the property $(N_i)$ for some integer $i\geq 1$, if
$F_j=\oplus S(-j-1)$ for all $j\leq i$ (or equivalently in terms
of graded Betti numbers, $b_{i,l}(C)=0$ for all $l\geq 2$). Using
the computation of $b_{j,l}(C)$ in terms of Koszul cohomology,
there is a well-known cohomological interpretation of property
$(N_i)$: we denote by $M=\Omega_{\PP^r}(1)$ and $M_L=M\otimes
\OO_C$, hence we have an exact sequence
$$ 0\rightarrow M_L\rightarrow H^0(L)\otimes \OO_C\rightarrow L\rightarrow
0.$$ By taking exterior powers, for each $i\geq 0$ we obtain the
exact sequence:
\begin{equation}
\label{wedge} 0\rightarrow \wedge^{i+1} M_L\rightarrow \wedge
^{i+i} H^0(L)\otimes \OO_C \rightarrow \wedge ^{i} M_L\otimes
L\rightarrow 0. \end{equation}
 If $L$ is a normally generated line
bundle on $C$, then $(C,L)$ satisfies property $(N_i)$ if and only
if for all $j\geq 1$, the natural map
$$u_{i, j}:\wedge^{i+1} H^0(L)\otimes H^0(L^{\otimes j })\rightarrow H^0(\wedge ^i
M_L\otimes L^{\otimes (j+1)})$$ obtained by tensoring the sequence
(\ref{wedge}) and taking global sections, is surjective (cf. e.g.
\cite{GL}, Lemma 1.10).

We will be interested in the vector bundle $M_L$ in the case when
$C$ is a curve of genus $g=2k-2$ and $L=K_C\otimes A^{\vee}$ is
 residual to a base point free pencil $A\in W^1_k(C)$. We start by 
establishing a more general technical
result used throughout the paper:
\begin{proposition}\label{vanishing}
Let $C$ be a curve of genus $g$ and $L$ a globally generated
$\mathfrak g^r_d$ on $C$. If $p\geq 1$ is an integer such that
$g+\rm{min}$$\{r, p\}\leq d$, then for all $j\geq 2$ and $0\leq
i\leq p+2-j$, we have the vanishing $H^1(\wedge^i M_L \otimes
L^{\otimes j})=0$.
\end{proposition}
\begin{proof} We use a filtration argument due to Lazarsfeld \cite{L}. 
We choose general points $x_1,\ldots, x_{r-1}\in C$
such that $h^0\bigl(C, L\otimes \OO_C(-x_1-\cdots-x_{r-1})\bigr)=2$
and then there is an exact sequence
$$0\longrightarrow L^{\vee}(x_1+\cdots+x_{r-1})\longrightarrow
M_L\longrightarrow \oplus_{j=1}^{r-1} \OO_C(-x_j)\longrightarrow
0.$$ Taking exterior powers, we can write the exact sequence
$$0\rightarrow \wedge^{i-1}\bigl(\oplus_{j=1}^{r-1} \OO_C(-x_j)\bigr)\otimes
L^{\vee}(x_1+\cdots +x_{r-1})\longrightarrow \wedge^i
M_L\longrightarrow \wedge^i \bigl(\oplus_{j=1}^{r-1}
\OO_C(-x_j)\bigr)\longrightarrow 0.$$ In order to conclude that
$H^0(\wedge^i M_L\otimes L^{\otimes j})=0$, it suffices to show that
for all $1\leq i\leq \mbox{min}\{r-1, p\}$ we have the following:

\noindent
(1) $H^0(L^{\otimes j}\otimes \OO_C(-D_i))=0$ for each effective divisor $D_i\in C_i$ with
support in the set $\{x_1, \ldots, x_{r-1}\}$, and 

\noindent
(2) $H^0\bigl(L^{\otimes (j-1)}\otimes \OO_C(D_{r-i})\bigr)=0$, for any effective divisor $D_{r-i}\in C_{r-i}$
with support contained in $\{x_1, \ldots, x_{r-1}\}$. 

\noindent
Statement $(2)$ follows because of degree reasons. The only way $(1)$ could fail is when $j=2$ and $h^0(K_C\otimes L^{\vee})\geq r-i+1$, which by the Riemann-Roch theorem implies that $g+i\geq d+1$ and this contradicts our numerical assumption.

\end{proof}

\begin{remark}
 
(1) For $p\geq r$, the condition $g+\rm{min}$$\{p,
r\}\leq d$ is equivalent with $H^1(L)=0.$
\newline  (2) Since we work in characteristic $0$, the same filtration argument from the proof of Proposition \ref{vanishing} can be used to show that $M_L$ is a stable vector bundle on $C$.
\end{remark}

We fix $g=2k-2$ where $k\geq 3$. We denote by $\cM_g^0$ the open
substack of $\cM_g$ consisting of curves $C$ of genus $g$ which
carry no $\mathfrak g^1_{k-1}$'s and which have no automorphisms
(clearly the complement of $\cM_g^0$ inside $\cM_g$ has
codimension $\geq 2$). Over $\cM_g^0$ we have a universal curve
$\pi:\cC\rightarrow \cM_g^0$ and we also consider the Hurwitz
stack $\sigma:\mathfrak G^1_k\rightarrow \cM_g^0$ parametrising
pairs $(C,A)$, where $[C]\in \cM_g^0$ and $A\in W^1_k(C)$ is a
(necessarily base point free) $\mathfrak g^1_k$ on $C$. It is a
classical result that $\mathfrak G^1_k$ is a smooth irreducible
stack which is finite over $\cM_{g}^0$. The isomorphism
$\mathfrak G^1_k\ni (C, A)\mapsto (C, K_C\otimes A^{\vee})\in
\mathfrak G^{k-2}_{3k-6}$ will be used throughout the paper.

Let us consider  now an element $(C,A) \in \mathfrak G^1_k$ and
set $L:=K_C\otimes A^{\vee}$. The genericity assumption on $C$
implies that $L$ is very ample and gives an embedding $C\subset
\PP^{k-2}$ of degree $3k-6$. We denote by $\I_C$ the ideal sheaf
of $C$ in this embedding.

\begin{proposition}\label{dege}
Fix integers $g=2k-2$ with $k\geq 3$ and $i\geq 0$ with $2k\geq
3i+7$.  If $(C,A)\in \mathfrak G^1_k$ and $L=K_C\otimes A^{\vee}$,
then $(C,L)$ fails to satisfy property $(N_i)$ if and only if
$$h^0(\PP^{k-2}, \wedge^i M_{\PP^{k-2}}\otimes \I_C(2))\geq
{k-2\choose i+1}(k-3i-6)+2{k-1\choose i+1}-{k-1\choose i+2}+1.$$
\end{proposition}

\begin{proof}
We use the following diagram of exact sequences for the embedding
$C\stackrel{|L|}\hookrightarrow \PP^{k-2}$:

$$\begin{array}{ccccc}
0 & \; & 0 & \;&  0\\
\rmapdown{} & \; & \rmapdown{} &\; & \rmapdown{}\\
 \wedge^{i+1}M_{\PP^{k-2}}\otimes\I_C(1) & \longrightarrow & \wedge^{i+1}H^0(\OO_{\PP^{k-2}}(1))\otimes \I_C(1) & \longrightarrow & \wedge^i M_{\PP^{k-2}}\otimes \I_C(2)  \\
 \rmapdown{} & \; & \rmapdown{} & \; & \rmapdown{}  \\
  \wedge^{i+1} M_{\PP^{k-2}}(1) & \longrightarrow  & \wedge^{i+1} H^0(\OO_{\PP^{k-2}}(1))\otimes \OO_{\PP^{k-2}}(1) & \longrightarrow & \wedge^i M_{\PP^{k-2}}(2) \\
 \rmapdown{a} & \; &\rmapdown{b} & \; &\rmapdown{c} \\
 \wedge^{i+1}M_L\otimes L &  \longrightarrow & \wedge^{i+1} H^0(L)\otimes L & \longrightarrow &
\wedge^i M_L\otimes L^2  \\
\rmapdown{} &\; & \rmapdown{} & \; & \rmapdown{}\\
0 & \;&  0 & \;&  0\\
\end{array}$$
For simplicity we denote by the same symbol the map on global sections corresponding to a sheaf morphism and the morphism itself. Since $H^1(\PP^{k-2}, \wedge^{i+1}M_{\PP^{k-2}}(1))=0$ (because of Bott vanishing, see \cite{OSS}, pg. 8), we can apply the snake lemma to the
diagram obtained by taking global sections in the two lower rows.
The map induced by $b$ is an isomorphism, which gives an
isomorphism
$$H^0(\PP^{k-2}, \wedge^i M_{\PP^{k-2}}\otimes \I_C(2))\cong H^1(\PP^{k-2}, \wedge^{i+1} M_{\PP^{k-2}}\otimes \I_C(1)).$$
From Proposition \ref{vanishing} we have that $H^1(\wedge^i
M_L\otimes L^2)=0$, and thus $C$
satisfies $(N_i)$ if and only if the map
$$H^1(\wedge^{i+1}M_L\otimes L)\rightarrow \wedge^{i+1}H^0(L)\otimes
H^1(L)$$ is an isomorphism, or equivalently $h^1(\wedge^{i+1}
M_L\otimes L)\leq h^1(\wedge^{i+1}H^0(L)\otimes L)=2{k-1\choose
i+1}$. From the Riemann-Roch theorem, this is the same thing as having
$$\mbox{dim}(\mbox{Coker}(a))\leq {k-2\choose i+1}(k-3i-6) +2{k-1\choose i+1}-{k-1\choose
i+2}.$$ But $\mbox{Coker}(a)\cong H^0(\wedge^i M_{\PP^{k-2}}\otimes
\I_C(2))$, so failure of $(N_i)$ is equivalent to
 the map $c$ having a kernel of dimension at least
 $$ {k-2\choose i+1}(k-3i-6)+2{k-1\choose i+1}-{k-1\choose i+2}+1=$$
$$=
 {k-2\choose i+1}(k-3i-6)\bigl(1-\frac{1}{(k-i-2)(i+2)}\bigr)+1.$$
\end{proof}

\begin{remark}\label{topstr}
 An immediate consequence of this proof is that when
$k-3i-6<0$ (or equivalently $g\leq 6i+8$), the pair $(C,A)$ always
fails property $(N_i)$. When $k=3i+6$ we have that $(C,A)$ fails
$(N_i)$ if and only if $h^0(\PP^{k-2}, \wedge^i
M_{\PP^{k-2}}\otimes \mathcal{I}_C(2))\geq 1$, which we expect to
be a  divisorial condition on $\mathfrak G^1_k$.
\end{remark}
This proposition allows us to define a determinantal substack of
$\mathfrak G^1_k$ consisting of pairs $(C,A)$ for which $L$ fails
to satisfy property $(N_i)$. On the fibre product
$\cC\times_{\cM_g^0}\mathfrak G^1_k$ there is a universal Poincar\'e
bundle $\L$ whose existence is guaranteed by the universal
property of $\mathfrak G^1_k$. If $p_1:\cC\times_{\cM_g^0}
\mathfrak G^1_k \rightarrow \cC$ and $p_2:\cC\times_{\cM_g^0}
\mathfrak G^1_k\rightarrow \mathfrak G^1_k$ are the natural
projections, then $\mathcal{E}:=p_{2 *}(p_1^*\omega_{\pi}\otimes
\L^{\vee})$ is a vector bundle of rank $k-1$ with fibre
$\E(C,A)=H^0(C, K_C\otimes A^{\vee})$ over each point $(C,A)\in
\mathfrak G^1_k$. We have a tautological embedding of the pullback
of the universal curve $\cC\times_{\cM_g^0}\mathfrak G^1_k$ into
the projective bundle $u:\PP(\E)\rightarrow \mathfrak G^1_k$ and
we denote by $\mathcal{J}\subset \OO_{\PP(\E)}$ the ideal sheaf of
the image. Next, we define the vector bundle $\cN$  on $\PP(\E)$
by the sequence
$$0\longrightarrow \cN\rightarrow u^*(\E)\longrightarrow \OO_{\PP(\E)}(1)\longrightarrow 0,$$
and we further introduce two vector bundles $\cA$ and $\cB$ over
$\mathfrak G^1_k$ by setting
$$\cA:=u_*\bigl(\wedge^i \cN\otimes \OO_{\PP(E)}(2)\bigr), \mbox{ and }
\cB:=u_*\bigl(\wedge ^i \cN\otimes \OO_{\cC\times _{\cM_g^0}
\mathfrak G^1_k}(2)\bigr).$$ If $C\subset \PP^{k-2}$ is the
embedding given by $L=K_C\otimes A^{\vee}$, then
$$\cA(C,A)=H^0(\PP^{k-2}, \wedge^i M_{\PP^{k-2}}(2)) \mbox{
}\mbox{ and }\cB(C,A)=H^0(C, \wedge^i M_L\otimes L^{\otimes 2})$$ and there
is a natural vector bundle morphism $\phi:\cA\rightarrow \cB$.
From Grauert's Theorem we see that both $\cA$ and $\cB$ are vector
bundles over $\mathfrak G^1_k$ (for $\cB$ use that $H^1(C,
\wedge^i M_L\otimes L^2)=0$ because of Proposition \ref{vanishing}). Moreover
from Bott's Theorem (see again \cite{OSS}) and the Riemann-Roch Theorem respectively, we obtain that
$$\mbox{rank}(\cA)=(i+1){k \choose i+2}\mbox{  and  }\mbox{
rank}(\cB)={k-2\choose i}(4k-9-3i).$$

Then Proposition \ref{dege} can be restated as follows:
\begin{proposition}\label{ni}
The cycle $$\mathcal{U}_{g, i}:=\{(C,A)\in \mathfrak G^1_k
:(C,K_C\otimes A^{\vee}) \mbox{ fails property }(N_i)\},$$ can be
realized as the degeneracy locus of rank $(i-1){k-1\choose
i+1}+(i+2){k-1\choose i+2}-{k-2\choose i+1}(k-3i-6)-1$ of the
vector bundle map $\phi:\cA\rightarrow \cB$ over $\mathfrak
G^1_k$.
\end{proposition}

In this way we obtain a stratification of $\cM_g^0$ with strata
$\mathcal{Z}_{g,i}:=\sigma(\mathcal{U}_{g, i})$ consisting of
those curves $C$ for which there exists a pencil $A\in W^1_{k}(C)$
such that $(C, K_C\otimes A^{\vee})$ fails property $(N_i)$. Note
that when $g=6i+10$ then $\mbox{rank}(\cA)=\mbox{rank}(\cB)$ and
$\mathcal{U}_{g,i}$ is simply the degeneracy locus of $\phi$ and
we expect $\mathcal{Z}_{g, i}$ to be a divisor on $\cM_g^0$.

 It is of course very natural to compare
this newly defined stratification
$$\mathcal{Z}_{g, 0}\subset \mathcal{Z}_{g, 1}\subset \ldots\subset
\mathcal{Z}_{g, i}\subset\ldots\subset \cM_g^0$$ with the more
classical stratification of $\cM_g$ defined in terms of syzygies
of $K_C$. If for each $i\geq 0$ we set
$$\cN_{g, i}:=\{[C]\in \cM_g:(C, K_C) \mbox{ fails property }(N_i)\},$$
then Green's Conjecture for generic curves of fixed gonality (cf.
\cite{V3}, \cite{V2}) can be read as saying that for all $0\leq
i\leq (g-2)/2$, the locus $\cN_{g, i}$ coincides with the
$(i+2)$-gonal locus $\cM_{g, i+2}^1$ of curves with a $\mathfrak
g^1_{i+2}$, while $\cN_{g, i}=\cM_g$ for all $i\geq (g-1)/2$. One
of the morals of this paper is that the stratification
$\mathcal{Z}_{g, i}$ is very different in nature from the one
given by gonality: whereas $\cN_{g, 0}$ is the locus of
hyperelliptic curves of genus $g$, the smallest stratum
$\mathcal{Z}_{g, 0}$ is intimately related to the locus of curves
lying on  $K3$ surfaces which is well-known to be transversal to
any Brill-Noether locus (a curve of genus $g$ lying on a general
$K3$ surface satisfies the Brill-Noether Theorem, cf. \cite{L}).

  We have already seen
 that $\mathcal{Z}_{g,i}=\cM_g^0$ when $g\leq 6i+8$ (cf. Remark \ref{topstr}).
  A
result due to Mukai and Voisin (cf. \cite{V1}, Proposition 3.2)
states that if $C$ is a
 smooth curve of genus $g=2k-2$ sitting on a $K3$ surface $S$ with
 $\mbox{Pic}(S)=\mathbb Z\ [C]$, then for all $A\in W^1_k(C)$ the
 multiplication map $\mbox{Sym}^2 H^0(K_C\otimes
 A^{\vee})\rightarrow H^0((K_C \otimes A^{\vee})^{\otimes 2})$ is not
 surjective, in other words, $[C]\in \mathcal{Z}_{g,0}$. In a
 forthcoming paper we will show that the converse also holds, that
 is, the closure in $\cM_g$ of the smallest stratum $\mathcal{Z}_{g,0}$
coincides with the locus
$$\mathcal{K}_g:=\{[C]\in \cM_g: C \mbox{ lies on  a }K3 \mbox{ surface}\}.$$
The possibility of such an equality of cycles has already been
raised in Voisin's paper (cf. \cite{V1}, Remarques 4.13). Its main
appeal lies in the fact that it gives an intrinsic
characterization of a curve lying on a $K3$ surface which makes no
reference to the $K3$ surface itself! We now make the following:
\begin{conjecture}\label{stratum}
 For an even genus $g\geq 6i+10$, the stratum $\mathcal{Z}_{g,i}$ is a proper
subvariety of $\cM_g^0$. In particular, when $g=6i+10$, the
stratum $\mathcal{Z}_{g,i}$ is a divisor on $\cM_g^0$.
\end{conjecture}

For $g\leq 8$ the conjecture is trivially true because
$\mathcal{K}_g=\cM_g$, that is, every curve of (even) genus $g\leq
8$ sits on a $K3$ surface. The first interesting case is $g=10$
when the conjecture holds and the identification
$\K_{10}=\mathcal{Z}_{10, 0}$ is part of a more general picture
(see \cite{FP} for details). We have checked Conjecture
\ref{stratum} for all genera $g\leq 24$ sometimes using the
program Macaulay. We will describe this in detail in the most
interesting cases, $g=16$ and $g=22$, when the strata
$\mathcal{Z}_{16, 1}$ and $\mathcal{Z}_{22, 2}$ are divisors on
$\cM_{16}$ and $\cM_{22}$ respectively.
\begin{theorem}\label{gen16}
Let $C$ be a general curve of genus $16$. Then every linear series
$A\in W^1_9(C)$ gives an embedding $C\stackrel{|K_C\otimes
A^{\vee}|}\longrightarrow \PP^7$ of degree $21$ such that the
ideal of $C$ is generated by quadrics (that is, it satisfies
property $(N_1)$). In particular $\mathcal{Z}_{16, 1}$ is a
divisor on $\cM_{16}$.
\end{theorem}
\begin{proof} From the irreducibility of the Hurwitz scheme of
coverings of $\PP^1$ it follows that the variety $\mathfrak
G^1_9=\mathfrak G^1_{16,9}$ parametrizing pairs $(C,A)$ with
$[C]\in \cM_{16}$ and $A\in W^1_9(C)$ is irreducible. To prove
that $\mathcal{Z}_{16, 1}$ is a divisor it is enough to exhibit an
element $(C,A)\in \mathfrak G^1_9$ with $K_C\otimes A^{\vee}$ very
ample, such that the embedded curve $C\stackrel{|K_C\otimes
A^{\vee}|} \hookrightarrow \PP^7$ is cut out by quadrics.

We consider $13$ general points in $\PP^2$ denoted by $p_1, p_2,
q_1, \ldots, q_7$ and $r_1,\ldots, r_4$ respectively, and define the
linear system
$$H\equiv 8h-3(E_{p_1}+E_{p_2})-2\sum_{i=1}^7 E_{q_i}-\sum_{j=1}^4
E_{r_j}$$ on the blow-up $S=\mbox{Bl}_{13}(\PP^2)$. Here $h$ denotes the
pullback of the line class from $\PP^2$. Using the program
Macaulay it is easy to check that $S\stackrel{|H|}\hookrightarrow
\PP^7$ is an embedding and the graded Betti diagram of $S$ is

\[
\begin{matrix}
 1 & -  & -  & -   \\
 - & 7 & - & -  \\
 - & -  & 35  & 56&
\end{matrix}
\]

Thus $S$ is cut out by quadrics.  To carry out this calculation we
chose the $13$ points in $\PP^2$ randomly using the Hilbert-Burch
theorem so that they satisfy the Minimal Resolution Conjecture
(see \cite{SchT} for details on how to pick random points in
$\PP^2$ using Macaulay). Next we consider a curve $C\subset S$ in
the linear system
 \begin{equation} \label{c16}
C\equiv 14h-5(E_{p_1}+E_{p_2})-4\sum_{i=1}^6 E_{q_i}-3E_{q_7}-
2\sum_{j=1}^3 E_{r_j}-E_{r_4}. \end{equation} By using Macaulay we
pick $C$ randomly in its linear system and then check that $C$ is
smooth, $g(C)=16$ and $\mbox{deg}(C)=21$. To show that $C$ is cut
out by quadrics one can compute directly the Betti diagram
of $C$. Otherwise, since $S$ is cut out by quadrics, to
conclude the same thing about $C$, it suffices to show that the
map
$$m:H^0(S, H)\otimes H^0(S, 2H-C)\rightarrow H^0(S, 3H-C)$$ is
surjective (or equivalently injective). Since $h^0(S,2H-C)=2$,
from the base point free pencil trick we get that
$\mbox{Ker}(m)=H^0(S, C-H)=0$, because $C-H\equiv
6h-2E_{p_1}-2E_{p_2}-2\sum_{i=1}^6 E_{q_i}-E_{q_7}-\sum_{j=1}^3
E_{r_j}$ is clearly not effective for a general choice of the $13$
points.
\end{proof}

\begin{remark}
A possible objection to this proof could be that while we defined
the Hurwitz space $\mathfrak G^1_9\rightarrow \cM^0_{16}$ over the
open part of $\cM_{16}$ consisting of curves carrying no
$\mathfrak g^1_8$'s, the curve $C$ we constructed in the proof of
Theorem \ref{gen16} might lie outside $\cM^0_{16}$. However, this
does not affect our conclusions because for every $g$ and
$l<(g+3)/2$, the Hurwitz space of pairs $(C, A)$ with $[C]\in
\cM_g$ an $A\in W^1_l(C)$, is an \emph{irreducible} variety, that
is, there can be no component of the Hurwitz space mapping to a
proper subvariety of the $l$-gonal locus $\cM_{g, l}^1$ (This
statement can be easily proved with the methods from \cite{AC1}).
\end{remark}
\begin{remark}\label{cusp16}
For later use we are going to record a slight generalization of
Theorem \ref{gen16}. The conclusion of the theorem holds for the
image in $\PP^7$ of any plane curve $C$ that is sufficiently
general in the linear system given by $(\ref{c16})$. In particular
one can choose $C$ to have a cusp at a general extra point $y\in
\PP^2$, in which case we obtain that a general one-cuspidal curve
$C\subset \PP^7$ of arithmetic genus $16$ and degree $21$ is cut
out by quadrics (this again is quite easily checked with
Macaulay). Since the Hilbert scheme of these curves is irreducible
(use again the Hurwitz scheme), it follows that for a
\emph{general} $(C, y)\in \cM_{15, 1}$ and for an \emph{arbitrary}
linear system $L\in W^7_{21}(C)$ having a cusp at $y$, the ideal
of the cuspidal image curve $C\stackrel{|L|}\rightarrow \PP^7$ is
generated by quadrics.
\end{remark}

In a somewhat similar manner we are going to show that for $g=22$
the locus $\mathcal{Z}_{22, 2}$ is a divisor.

\begin{theorem}\label{gen22}
Let $C$ be a general curve of genus $22$. Then for every linear
series $A\in W^1_{12}(C)$ we have that $|K_C\otimes A^{\vee}|$ is
very ample and the resulting embedding $A\stackrel{|K_C\otimes
A^{\vee}|}\hookrightarrow \PP^{10}$ satisfies property $(N_2)$. In
particular, the locus of curves $C$ of genus $22$ carrying a linear
series $\mathfrak g^{10}_{30}=K_C(-\mathfrak g^1_{12})$ which fails
to satisfy property $(N_2)$, is a divisor on $\cM_{22}$.
\end{theorem}
\begin{proof} Ideally we would like to construct a rational surface
$S\subset \PP^{10}$ which satisfies property $(N_2)$ and then
consider a suitable curve $C\subset S$ with $g(C)=22$ and
$\mbox{deg}(C)=30$. For numerical reasons this will turn out not
to be possible but we will be able to construct a smooth curve
$C'\subset S$ with $g(C')=20$ and $\mbox{deg}(C')=28$ which
satisfies $(N_2)$. The desired curve $C$ will be a smoothing in
$\PP^{10}$ of $C\cup L\cup L'$, where $L$ and $L'$ are general
chords of $C'$.

We start by choosing random points $p, p_1,\ldots, p_7, q_1, \ldots,
q_5\in \PP^2$. Then the linear system
$$H=11h-4E_p-3\sum_{i=1}^7 E_{p_i}-2\sum_{j=1}^5 E_{q_j}$$
yields an embedding $S\subset \PP^{10}$ of $\mbox{Bl}_{13}(\PP^2)$. Using
Macaulay we get that the upper left corner of the Betti diagram of
$S$ is
\[
\begin{matrix}
 1 & -  & -  & -  & - \\
 - & 29 & 98 & 72 & - \\
 - & -  & -  & 264& -
\end{matrix}
\]
that is,  $S$ satisfies property $(N_2)$, hence $H^1(S, \wedge^3
M_S(1))=0$. If $\Gamma\equiv H$ is a general hyperplane section of
$S$, then $g(\Gamma)=13$, $\mbox{deg}(\Gamma)=22$ and an argument
very similar to the one in Proposition \ref{vanishing} shows that the
vector bundle $M_{\Gamma}$ is stable. In particular $H^1(\Gamma,
\wedge^2 M_{\Gamma}(2))=0$. Since $M_{S | \Gamma}\cong
M_{\Gamma}\oplus \OO_{\Gamma}$ one can write the exact sequence
$$0\longrightarrow \wedge^{i}M_S(1)\longrightarrow \wedge^i M_S(2)\longrightarrow \wedge^i
M_{\Gamma}(2)\oplus \wedge^{i-1} M_{\Gamma}(2)\longrightarrow 0,$$
which leads to the vanishing $H^2(S, \wedge^3 M_S(1))=0$. Suppose
now that $C\subset S$ is a smooth non-degenerate curve. We have a
commutative diagram:
$$\begin{array}{cccccc} \; & \; & H^1\bigl(C, \wedge^3
M_C(1)\bigr) & \longrightarrow & \wedge^3
H^0\bigl(\OO_C(1)\bigr)\otimes H^1(\OO_C(1))  \\
\; & \; & \rmapdown{=} & \; & \rmapdown{=}  \\
H^1\bigl(\wedge^2 M_S(2H-C)\bigr) & \hookrightarrow  & H^2\bigl(
\wedge^3 M_S(H-C)\bigr) &\rightarrow & \wedge^3 H^0(
\OO_S(1))\otimes
H^2(H-C) ,\\
\end{array}$$
where the left vertical map is isomorphic because $H^1(S, \wedge^3
M_S(1))=H^2(S, \wedge^3 M_S(1))=0$. To conclude that $C$ satisfies
$(N_2)$, it suffices to show that $H^1\bigl(S, \wedge^2
M_S(2H-C)\bigr)=0$, or equivalently, that the map $f:\wedge^2
H^0(\OO_S(1))\otimes H^0(2H-C)\rightarrow H^0(M_S(3H-C))$ obtained
from the Koszul complex, is surjective. If $\bigl(\mbox{deg}(C),
g(C)\bigr)=(30, 22)$, then one can check easily that $h^0(S,
2H-C)=3$ and $h^0(S, M_S(3H-C))=207$, hence $f$ cannot be
surjective for dimensional reasons. The closest we can get to
these numerical invariants is when $(\mbox{deg}(C),g(C))=(28,20)$
and this is the type of curve on $S$ we will be looking for. Take
$$C'\equiv 17h-6E_p-5\sum_{i=}^7 E_{p_i}-3\sum_{j=1}^5 E_{q_j},$$
and consider the embedding $C'\subset \PP^{10}$ given by $|H|$.
Using Macaulay we check that $C'$ is a smooth curve of genus $20$
and degree $28$ with graded Betti diagram
\[
\begin{matrix}
 1 & - & - & - \\
 - & 27 & 80  & - \\
 - & -  & - & 432
\end{matrix}
\]
Hence $C'$ satisfies $(N_2)$. Now choose two general chords $L$ and
$L'$ of $C$ and define $C:=C'\cup L\cup L'$. Clearly $g(C)=22,
\mbox{deg}(C)=30$ and one last check with Macaulay shows that
$b_{2j}(C)=0$ for $j\geq 2$, that is, $C$ satisfies $(N_2)$.

\end{proof}

\begin{remark}\label{cusp22}

Just like in the case $g=16$ we have a slight variation of the last
Theorem. The same proof shows that for a \emph{general} $(C, y)\in
\cM_{21, 1}$ and for an \emph{arbitrary} linear series $L\in
W^{10}_{30}(C)$ that has a cusp at $y$, the one-cuspidal image curve
$C\stackrel{|L|}\rightarrow \PP^{10}$ satisfies property $(N_2)$.
\end{remark}

\section{Intersection theory calculations on $\mm_g$}

Recall that we have realized $\mathcal{Z}_{6i+10, i}$ as the image
of the degeneracy locus $\mathcal{U}_{6i+10, i}$ of a vector bundle
morphism $\phi:\cA\rightarrow \cB$ over $\mathfrak G^1_k$. To
compute the class of the compactification $\overline{\mathcal{Z}}_{6i+10, i}$  we are
going to extend the determinantal structure of
$\mathcal{Z}_{6i+10, i}$ over the boundary divisors in
$\mm_{6i+10}$.

We set $g:=6i+10, k:=3i+6$ and denote by
$\widetilde{\cM}_g:=\cM_g^0\cup \bigl(\cup_{j=0}^{3i+5}
\Delta_j^0\bigr)$ the locally closed substack of $\mm_g$ defined as
the union  of the locus $\cM_{g}^0$ of smooth curves carrying no
linear systems $\mathfrak g^1_{k-1}$  to which we add the open
subsets $\Delta_j^0\subset \Delta_j$ for $1\leq j\leq 3i+5$
consisting of $1$-nodal genus $g$ curves $C\cup_y D$, with $[C]\in
\cM_{g-j}$ and $[D, y]\in \cM_{j, 1}$ being Brill-Noether general
curves, and the locus $\Delta_0^0\subset \Delta_0$ containing
$1$-nodal irreducible genus $g$ curves $C'=C/q\sim y$, where $[C,
q]\in \cM_{g-1}$ is a Brill-Noether general pointed curve and
$y\in C$, together with their degenerations consisting of unions
of a smooth genus $g-1$ curve and a nodal rational curve. One can
then extend the finite covering $\sigma:\mathfrak
G^{k-2}_{3k-6}\rightarrow \cM_g^0$ to a proper, generically finite
map
$$\sigma: \widetilde{\mathfrak G}^{k-2}_{3k-6} \rightarrow \widetilde{\cM}_g$$ by
letting $\widetilde{\mathfrak G}^{k-2}_{3k-6}$ be the variety of
limit $\mathfrak g^{k-2}_{3k-6}$'s on the treelike curves from
$\widetilde{\cM}_g$ (see \cite{EH1}, Theorem 3.4 for the
construction of the space of limit linear series).

We will be
interested in intersecting the divisors
$\overline{\mathcal{Z}}_{6i+10, i}$ on $\mm_{g}$ with test
curves in the boundary of $\mm_{6i+10}$ which are defined as
follows: we fix a Brill-Noether general curve $C$ of genus
$2k-3=6i+9$, a general point $q\in C$ and a general elliptic curve
$E$. We define two $1$-parameter families
\begin{equation}\label{testcurves}
C^0:=\{\frac{C}{y\sim q}: y\in C\}\subset \Delta_0 \subset \mm_{2k-2}
\mbox{ and }C^1:=\{C\cup _y E: y\in C\}\subset \Delta_1\subset
\mm_{2k-2}.
\end{equation}
 It is well-known that
these families intersect the generators of
$\mbox{Pic}(\mm_{2k-2})$ as follows:
$$
C^0\cdot \lambda=0,\ C^0\cdot \delta_0=-(4k-6), \ C^0\cdot
\delta_1=1 \mbox{ and } C^0\cdot \delta_a=0\mbox{ for }a\geq 2,
\mbox{ and}$$
$$C^1\cdot \lambda=0, \ C^1\cdot \delta_0=0, \ C^1\cdot
\delta_1=-(4k-8), \ C^1\cdot \delta_a=0 \mbox{ for }a\geq 2.$$
Next, we fix an integer $2\leq j\leq k-1$, a general curve $C$ of
genus $2k-2-j$ and a general curve pointed curve $(D, y)$ of
genus $j$. We define the $1$-parameter family $C^j:=\{C\cup_y D:
y\in C\}\subset \Delta_j\subset \mm_{2k-2}$. We have that
$$C^j\cdot \lambda=0, \ C^j\cdot \delta_a=0 \mbox{ for }a\neq j
\mbox{ and } C^j\cdot \delta_j=-(4k-6-2j).$$
 During our
calculations we will often need the following enumerative result (cf.
\cite{HM}, pg. 71):
\begin{proposition}\label{enum}
Let $C$ be a general curve of genus $j$, $y\in C$ a general fixed
point and an integer $0\leq \alpha \leq j/2.$
\begin{enumerate} \item
The number of pencils $A \in W^1_{j-\alpha+1}(Y)$ satisfying
$h^0(Y, A(-(j+1-2\alpha)y))\geq 1$ is equal to $$a(j,
\alpha)=\frac{j+1-2\alpha}{j+1-\alpha}{j \choose \alpha}.$$ \item
The number of pencils $A\in W^1_{j-\alpha}(C)$ satisfying $h^0(Y,
A(-(j-2\alpha)q))\geq 1$ for some unspecified point $q\in Y$ is
equal to
$$b(j, \alpha)=(j-2\alpha-1)(j-2\alpha)(j-2\alpha+1){j \choose
\alpha}.$$
\end{enumerate}
\end{proposition}

\begin{definition}
(1) For a fixed general curve $C$ of genus $2k-3$ we denote by
$\delta (k, 1)$ the (finite) number of linear series $L\in
W^{k-2}_{3k-6}(C)$ for which  there exists a point $y\in C$ with
$h^0(L(-2y))\geq k-2$ and $h^0(L(-ky))\geq 1$.
\newline \noindent (2) We fix $2\leq j\leq k-1,\ 0\leq \alpha\leq j/2$ and a general
curve $C$ of genus $2k-2-j$. We denote by $\delta(k, j, \alpha)$
the number of linear series $L_C \in W^{k-2}_{3k-j-4}(C)$ for
which there exists a point $y\in C$ such that
$h^0\bigl(L(-(\alpha+1)y)\bigr)\geq k-1-\alpha$, $h^0\bigl(L(-(j-\alpha+2)y)\bigr)\geq
k-1+\alpha-j$ and $h^0\bigl(L(-(k+1)y)\bigr)\geq 1$. \end{definition}

 The numbers $\delta(k, 1)$ and $\delta(k, j, \alpha)$ can be computed by degeneration.
 Since we do not need their actual values we skip this calculation.
 Next we describe $\widetilde{\mathfrak
G}^{k-2}_{3k-6}$ set theoretically. To set notation, if $X$ is a
treelike curve and $l$ is a limit $\mathfrak g^r_d$ on $X$, for a
component  $Y$ of $X$ we denote by $l_Y=(L_Y, V_Y\subset
H^0(L_Y))$ the $Y$-aspect of $l$. For a point $y\in Y$ we denote
by  $\{a^{l_Y}_s(C)\}_{s=0\ldots r}$ the \emph{vanishing
sequence} of $l$ at $y$ and by $\rho(l_Y, y):=\rho(g, r,
d)-\sum_{i=0}^r (a^{l_Y}_i(y)-i)$ the adjusted Brill-Noether
number with respect to $y$
\begin{proposition}\label{limitlin}
(1) Let $C_y^1=C\cup_y E$ be an element of $\Delta_1^0$. If $(l_C,
l_E)$ is a limit  $\mathfrak g^{k-2}_{3k-6}$ on $C_y^1$, then
$V_C=H^0(L_C)$ and $L_C\in W^{k-2}_{3k-6}(C)$ has a cusp at $y$,
that is, $h^0(C, L(-2y))=k-2$. If $y\in C$ is a general point, then
$l_E=\bigl(\OO_E((3k-6)y), (2k-5)y+|(k-1)y|\bigr)$, that is, $l_E$
is uniquely determined. If $y\in C$ is one of the finitely many
points for which there exists $L_C\in W^{k-2}_{3k-6}(C)$ such that
$\rho(L_C, y)=-1$, then $l_E(-(2k-6)y)$ is a $\mathfrak g_k^{k-2}$
with vanishing sequence at $y$ being $\geq (0, 2, 3, \ldots, k-2,
k)$.  Moreover, at the level of $1$-cycles we have the
identification $\sigma^*(C^1)\equiv X + \delta (k, 1)\ T$, where
$$X:=\{(y, D)\in C\times C_{k-2}:h^0(C, D+2y)\geq 2\}$$
and $T$ is the curve consisting of $\mathfrak g^{k-2}_k$'s on $E$
with vanishing $\geq (0, 2, \ldots, k-2, k)$ at the fixed point
$y\in E$.

\noindent (2) Let $C_y^0=C/y\sim q$ be an element of $\Delta_0^0$.
Then limit linear series of type $\mathfrak g^{k-2}_{3k-6}$ on
$C_y^0$ are in 1:1 correspondence with complete linear series $L$
on $C$ of type $\mathfrak g^{k-2}_{3k-6}$ satisfying the condition
$h^0(C, L\otimes \OO_C(-y-q))=h^0(C,L)-1.$ Thus there is an
isomorphism between the cycle $\sigma^*(C^0)$ of $\mathfrak
g^{k-2}_{3k-6}$ on all curves $C_y^0$ with $y\in C$ and the smooth
curve
$$Y:=\{(y,D)\in C\times C_{k-2}: h^0(C, D+y+q)\geq 2\}.$$

\end{proposition}

\begin{proof} We only prove (1) the remaining case being similar.
Suppose $l=(l_C, l_E)$ is a limit $\mathfrak g^{k-2}_{3k-6}$ on $C_y^1$. Using the additivity of the
Brill-Noether number (see e.g. \cite{EH1}, Lemma 3.6), we get that
$0=\rho(2k-2, k-2, 3k-6)\geq \rho(l_C,y)+\rho(l_E,y)$. Since the
vanishing sequence at $y$ of the $E$-aspect of $l$ is
$$a^{l_E}(y)\leq (2k-5,2k-4,\ldots,  3k-8, 3k-6),$$ it follows
that we also have the inequality $a^{l_C}(y)\geq (0, 2,\ldots,
k-1)$. If $y\in C$ is general, then $\rho(l_C, y)=0$ and
$a^{l_C}(y)=(0,2,\ldots, k-1)$. Moreover, $l_C$ corresponds to a
complete linear series $|L_C|$ on $C$ such that $h^1(C,L_C)=1$ and
$h^1(L_C\otimes \OO_C(-2y))=2$. If $y\in C$ is one of the points
such that $\rho(l_C, y)=-1$, then since $C$ is Brill-Noether
general we must have $a^{l_C}_0(y)=0$, which implies that
$a^{l_C}(y)=(0, 2, \ldots, k-2, k)$, hence $a^{l_E}(y)\geq (2k-6,
2k-4, \ldots, 3k-8, 3k-6)$, or equivalently, $l_E=(2k-6)y+|V_E|$,
where $V_E\subset H^0(\OO_E(ky))$ is a $\mathfrak g^{k-2}_k$ with
vanishing sequence $\geq (0, 2, 3, \ldots, k-2, k)$ at $y$. An
easy argument shows that the variety $Y$ of such $V_E$'s is
isomorphic to $\PP(\OO_{E}(ky)_{| 2y})$.
\end{proof}

\begin{proposition}\label{limitlin2}
 Let $C_y^j:=C \cup _y D$ be an element from
$\Delta_j^0$, where $2\leq j\leq k-1$, $g(C)=2k-2-j$ and $g(D)=j$.
 Suppose $(l_C, l_D)$ is a limit $\mathfrak g^{k-2}_{3k-6}$ on
$C_y^j$. If $y\in C$ is a general point, then $l_C$ has the
divisor $(j-2)y$ as base locus and there exists an integer $0\leq
\alpha \leq [j/2]$ such that $l_C(-(j-2)y)=|L_C|$, with $L_C \in
W^{k-2}_{3k-j-4}(C)$ satisfying the inequality
$$a^{L_C}(y)\geq (0, 1, \ldots, \alpha-1,
\alpha+1, \alpha+2, \ldots, j-\alpha, j-\alpha+2,j-\alpha+3, \ldots,
k-1, k).$$ Moreover $l_D$ is one of the $a(j, \alpha)$ linear
systems on $D$ of type $(2k-j-4)y+|L_D|$, with $L_D\in
W^{k-2}_{k-2+j}(D)$ such that $h^0(L_D(-(k+1-\alpha)y))=\alpha
\mbox{ and }h^0(L_D(-(k+\alpha-j)y))=j-\alpha.$

If for some integer  $0\leq \alpha \leq [j/2]$,  $y\in C$ is one of
the $b(2k-2-j, k-1+\alpha-j)$ points for which there exists  $L_C\in
\rm{Pic}$$^{3k-j-3}(C)$ such that
$h^0(L_C(-(\alpha+2)y))=k-1-\alpha$ and
$h^0(L_C((j-\alpha+2)y))=k-j+\alpha$, then there are two
possibilities:

\noindent (1) The aspect $l_C$ is of the form $(j-2)y+|L_C(-y)|$, whereas $l_D=(2k-j-4)y+|L_D|$, where
$L_D\in W^{k-2}_{k-2+j}(D)$ satisfies
$h^0\bigl(L_D(-(k+1-\alpha)y)\bigr)\geq \alpha \mbox{ and
 }h^0\bigl(L_D(-(k+1+\alpha-j)y)\bigr)\geq j-\alpha-1.$

\noindent (2) The aspect $l_C$ is of the form $(j-3)y+(L_C, V_C)$, where $V_C\subset H^0(L_C)$ is such that $H^0(L_C(-2y))\subset V_C$, and $l_D=(2k-j-4)y+|L_D|$ with $L_D\in W_{k-2+j}^{k-2}(D)$ and
$$a^{L_D}(y)=(0, 1, \ldots, k-1-j+\alpha, k+1-j+\alpha, \ldots, k-1-\alpha, k+1-\alpha, \ldots, k-1, k+1).$$

Finally, if for some $0\leq \alpha \leq [j/2]$, $y\in C$ is  one of the $\delta(k, j, \alpha)$ points for which there exists $L_C\in \mbox{Pic}^{3k-4-j}(C)$  such that $$a^{L_C}(y)=(0, 1,
\ldots, \alpha-1, \alpha+1, \ldots, j-\alpha, j-\alpha+2, \ldots,
k-2, k-1, k+1),$$ then $l_C=(j-2)y+|L_C|$ and $l_D=(2k-j-5)y+l_D^{'}$, where $l_D^{'}$ is a (non-complete)
$\mathfrak g^{k-2}_{k-1+j}$ on $D$ satisfying
$a^{l_D^{'}}(y)=(0, 2, \ldots, k-1+\alpha-j, k+1+\alpha-j, \ldots,
k-\alpha, k+2-\alpha,\ldots, k, k+1)$.
\end{proposition}
\begin{proof} Since by definition $[D, y]\in \cM_{j, 1}$ is Brill-Noether
general, we have that $\rho(l_D, y)\geq 0$. The locus of those
$(C, y)\in \cM_{2k-2-j, 1}$ carrying a $\mathfrak g^{k-2}_{3k-6}$
with $\rho(\mathfrak g^{k-2}_{3k-6}, y)\leq -2$ is of codimension
$\geq 2$, hence we must also have that $\rho(l_C, y)\geq -1$.
Using Proposition 1.2 from \cite{EH3} we have that
$a_{k-2}^{l_C}(y)\leq k-2+j$ and $a_0^{l_C}(y)\geq j-3$. There are
two cases to consider:

\noindent {\bf (a)} If $\rho(l_C, y)=0$, then $\rho(l_D, y)=0$.
Assume first that $a_0^{l_C}(y)=j-2$. Since we have that $0\leq
a_s^{l_C}(y)-s-(j-2)\leq 2$ for all $0\leq s\leq k-2$, there must
exist an integer $0\leq \alpha\leq [j/2]$, such that if  $L_C$
denotes the $\mathfrak g^{k-2}_{3k-4-j}$ obtained from $l_C$ by
removing $(j-2)y$, then the vanishing at $y$ is $a^{L_C}(y)=(0, 1,
\ldots, \alpha-1, \alpha+1, \alpha+2, \ldots, j-\alpha, j-\alpha+2,
\ldots, k)$. By compatibility, $l_D$ has $(2k-j-4)y$ in its base
locus and we write $l_D=(2k-j-4)y+|L_D|$, where $a^{L_D}(y)=(0, 1,
\ldots, k-2+j-\alpha, k+j-\alpha, \ldots, k-1-\alpha, k+1-\alpha,
\ldots, k)$. By Riemann-Roch, the number of such $L_D$'s equals the
number of $\mathfrak g^1_{j-\alpha+1}$'s on $D$ with a
$(j-2\alpha+1)$-fold point at $y$ which is computed by $a(j,
\alpha)$. If $a_0^{l_C}(y)=j-3$, we denote by $L_C\in
\mbox{Pic}^{3k-3-j}(C)$ the line bundle obtained from $l_C$ by
subtracting $(j-3)y$, and again we see that there must exist $0\leq
\alpha \leq [j/2]$ such that $h^0(L_C(-(j-\alpha+2)y))=k-j+\alpha$
and $h^0(L_C(-(\alpha+2)y))=k-1-\alpha$. By duality, $K_C\otimes
L_C^{\vee}\otimes \OO_C((j-\alpha+2)y)$ is one of the finitely many
$\mathfrak g^1_{k-1-\alpha}$'s on $C$ with a $(j-2\alpha)$-fold
point at $y$. This number equals $b(2k-2-j, k-1+\alpha-j)$ and this
explains case (2) of our statement.

\noindent {\bf (b)} If $\rho(l_C, y)=-1$, then $\rho(l_D, y)=1$.
An argument similar to that used above gives two possibilities
for $a^{l_C}(y)$ and the rest follows because $l_C$ and $l_D$ are
compatible at $y$.
\end{proof}
Retaining the notation from Proposition \ref{limitlin2}, for each
integer $0\leq \alpha \leq [j/2]$ we define the following
$1$-cycles:
\newline
$$X_{j, \alpha}:=\{(y, L_C)\in C\times \mbox{Pic}^{3k-j-4}(C)
:a^{L_C}_{\alpha}(y)\geq \alpha+1,\mbox{ }
a^{L_C}_{j-\alpha}(y)\geq j-\alpha+2\},$$
$$Y_{j, \alpha}^{'}:=\{(y,l)\in C\times G^{k-2}_{3k-j-3}(C):
a^{l}_{1}(y)\geq 2, \mbox{ }
a^{l}_{\alpha}(y)\geq \alpha+2, \mbox{ } a^{l}_{j-\alpha-1}(y)\geq j-\alpha+2\},$$
$$Y_{j,
\alpha}^{''}:=\{L_D\in \mbox{Pic}^{k-2+j}(D):
a^{L_D}_{k+\alpha-j}(y)\geq k+\alpha-j+1, \mbox{
}a^{L_D}_{k-\alpha-1}(y)\geq k-\alpha+1 \}$$ and $$Y_{j,
\alpha}^{'''}:=\{l\in G^{k-2}_{k-1+j}(D): a_1^{l}(y)\geq 2,
a_{k-1+\alpha-j}^l(y)\geq k+\alpha-j+1, a^l_{k-1-\alpha}(y)\geq
k-\alpha+2\}.$$
Via  Proposition \ref{limitlin2}
we view $X_{j, \alpha}, Y_{j, \alpha}^{'}$, $Y_{j,
\alpha}^{''}$ and $Y_{j, \alpha}^{'''}$ as $1$-cycles on
$\widetilde{\mathfrak G}^{k-2}_{3k-6}$. Note that $X_{j, \alpha}$ projects onto $C$ while $Y_{j, \alpha}^{'}$ is isomorphic to the union of $b(2k-2-j, k-1+\alpha-j)$ copies of $\PP^1$. We can then rephrase
Proposition \ref{limitlin2} as follows:

\begin{proposition}\label{limitlin3}
If $\sigma:\widetilde{\mathfrak G}^{k-2}_{3k-6}\rightarrow
\widetilde{\cM}_g$ is the natural projection  and $[C]\in
\cM_{2k-2-j}$ where $2\leq j\leq k-1$, then we have the following
numerical equivalence relation between $1$-cycles:

$$\sigma^*(C^j)\equiv \sum_{\alpha =0}^{[j/2]} \Bigl(a(j,
\alpha)\  X_{j, \alpha}+\delta(k, 2k-2-j, k+\alpha-j)\ Y_{j,
\alpha}^{'}+$$ $$+b(2k-2-j, k-1+\alpha-j)\ Y_{j, \alpha}^{''}
+\delta(j, \alpha)\ Y_{j, \alpha}^{'''}\Bigr).$$
\end{proposition}

To compute the class of the cycles introduced in Propositions
\ref{limitlin} and \ref{limitlin3}, we introduce some notation.
Let $C$ be a Brill-Noether general curve of genus $g$. We denote
by $C_d$ the $d$-th symmetric product of $C$, by $\cU\subset
C\times C_d$ the universal divisor, and by $\L$ the Poincare
bundle on $C\times \mbox{Pic}^d(C)$. We also introduce the
projections $\pi_1:C\times C_d\rightarrow C$ and $\pi_2:C\times
C_d\rightarrow C_d$ (and we use the same notation for the
projections from $C\times \mbox{Pic}^d(C)$ onto the factors). We
define the cohomology classes
$$\eta=\pi_1^*([point])\in H^2(C\times C_d) \mbox{ and }
x=\pi_2^*(D_{x_0}) \in H^2(C\times C_d),$$ where $x_0\in C$ is an
arbitrary point and $D_{x_0}:=\{D\in C_d:x_0\in D\}$ (the
definition of $x$ is obviously independent of the point $x_0$).
Finally, if $\delta_1,\ldots, \delta_{2g}\in H^1(C, \mathbb
Z)\cong H^1(C_d, \mathbb Z)$ is a symplectic basis, then we define
the class
$$\gamma:=-\sum_{\alpha=1}^g
\Bigl(\pi_1^*(\delta_{\alpha})\pi_2^*(\delta_{g+\alpha})-\pi_1^*(\delta_{g+\alpha})\pi_2^*(\delta_
{\alpha})\Bigr).$$ With these notations we have the formula (cf.
\cite{ACGH}, p. 338) $[\cU]\equiv d\eta+\gamma+x,$ corresponding
to the Hodge decomposition of $[\cU]$. We also record the
formulas $\gamma^3=\gamma \eta=0$, $\eta^2=0$ and $\gamma^2=-2\eta
\pi_2^*(\theta)$, where $\theta \in H^2(C_d, \mathbb Z)$ is the
pullback of the class of the theta divisor on $\mbox{Jac}(C)$.
Similarly, $c_1(\L)=d \eta+\gamma$, where $\eta, \gamma \in
H^2(C\times \mbox{Pic}^d(C))$ pull back to $\eta$, $\gamma \in
H^2(C\times C_d).$ When $d=k-2$ and $g=2k-3$, we set
$\cM:=\pi_1^*(K_C)\otimes \OO_{C\times C_{k-2}}(-\cU)$, hence
$c_1(\cM)=(3k-6)\eta-\gamma-x$.

We now compute the class of the curves $X$ and $Y$ defined in
Proposition \ref{limitlin}:

\begin{proposition}\label{xy}
Let $C$ be a Brill-Noether general curve of genus $g=2k-3$ and $q\in
C$ a general point.

\noindent(1) The class of the curve $X=\{(y,D)\in C\times
C_{k-2}:h^0(C,D+2y)\geq 2\}$ is given by
$$[X]\equiv \frac{\theta^{k-2}}{(k-2)!}-\frac{\theta^{k-3}x}{(k-3)!}-\frac{2\theta^{k-3}\gamma}{(k-3)!}+$$
$$
 +\frac{2 x\theta^{k-4}\gamma}{(k-4)!} -\frac{4(k+1)}{(k-4)!}
x\theta^{k-4}\eta + \frac{4k-2}{(k-3)!}\theta^{k-3}\eta .
$$

\noindent (2) The class of the curve $Y=\{(y,D)\in C\times
C_{k-2}:h^0(C,D+y+q)\geq 2\}$ is given by
$$[Y]\equiv \frac{\theta^{k-2}}{(k-2)!}-\frac{\theta^{k-3}x}{(k-3)!}-\frac{\theta^{k-3}\gamma}{(k-3)!}+$$
$$+
\frac{x\theta^{k-4}\gamma}{(k-4)!}-\frac{k+1}{(k-4)!}
x\theta^{k-4}\eta+\frac{k-1}{(k-3)!}\theta^{k-3}\eta.$$
\end{proposition}
\begin{proof}
We will realize both $X$ and $Y$ as pullbacks of degeneracy loci
and compute their classes using the Thom-Porteous formula. We
consider the map $\epsilon:C\times C_{k-2}\rightarrow C_k$ given
by $\epsilon(y,D):=2y+D$. It is easy to check that
$$\epsilon^*(x)=2\eta+x \mbox{ and
}\epsilon^*(\theta)=4g\eta+\theta-2\gamma.$$ Then
$X=\epsilon^*(C^1_k)$ and to compute the class of $C^1_k=\{E\in
C_k: h^0(E)\geq 2\}$ we introduce the rank $k$ vector bundle
$\mathcal{P}:=(\pi_2)_*(\pi_1^*K_C\otimes \OO_{\mathcal{U}})$ on
$C_k$ having total Chern class
$c_t(\mathcal{P})=(1-x)^{-k+4}e^{\frac{\theta}{1-x}}$ (cf.
\cite{ACGH}, pg. 240). There is a natural bundle map
$\phi^{\vee}:\epsilon^*(\mathcal{P}^{\vee})\rightarrow
H^0(K_C)^{\vee}\otimes \OO_{C\times C_{k-2}}$, and
$X=Z_{k-1}(\phi^{\vee})$. Then by Thom-Porteous, we can write
\newline
\noindent
$$[X]=\bigl[\frac{1}{c_t(\epsilon^*(\mathcal{P}^{\vee}))}\bigr]_{k-2}=
\epsilon^*\bigl(\frac{\theta^{k-2}}{(k-2)!}-\frac{x\theta^{k-3}}{(k-3)!}\bigr)=$$
$$
=\frac{(4g\eta+\theta-2\gamma)^{k-2}}{(k-2)!}-\frac{(4g\eta+\theta-2\gamma)^{k-3}(2\eta+x)}{(k-3)!},
$$ which quickly leads to the desired expression for $[X]$.

To compute the class of $Y$ we proceed in a similar manner: we
consider the map $\chi:C\times C_{k-2}\rightarrow C_k$ given by
$\chi(y, D):=y+q+D$. Then $\chi^*(x)=x+\eta$ and
$\chi^*(\theta)=g\eta+\theta-\gamma$. For each $(y, D)\in C\times
C_{k-2}$ we have a natural map $$H^0(K_{C
|y+q+D})^{\vee}\rightarrow H^0(K_C)^{\vee}$$ which globalizes to a
vector bundle map
$\psi^{\vee}:\chi^*(\mathcal{P}^{\vee})\rightarrow
H^0(K_C)^{\vee}\otimes \OO_{C\times C_{k-2}}$ and then it is clear
that $Y=Z_{k-1}(\psi^{\vee})$. Applying Thom-Porteous again, we
obtain
\newline
\noindent
$$[Y]=\bigl[\frac{1}{c_t(\chi^*(\mathcal{P}^{\vee}))}\bigr]_{k-2}=\chi^*\bigl(\frac{\theta^{k-2}}{(k-2)!}
-\frac{x\theta^{k-3}}{(k-3)!}\bigr)=$$
$$=\frac{(g\eta+\theta-\gamma)^{k-2}}{(k-2)!}-\frac{(g\eta+\theta-\gamma)^{k-3}(\eta+x)}{(k-3)!},$$
which after a straightforward calculation gives the class of $Y$.
\end{proof}

We also need the following intersection theoretic result:
\begin{lemma}\label{fj}
For $j\geq 1$ we denote by $\F_j:=(\pi_2)_*(\cM^{\otimes j})$ the
vector bundle on $C_{k-2}$ with fibre $\F_j(D)=H^0(C,jK_C-jD)$
over each $D\in C_{k-2}$. Then
$c_1(\F_j)=-j^2\theta-j(3j-2)(k-2)x$.
\end{lemma}
\begin{proof} We apply the Grothendieck-Riemann-Roch Theorem for the
map $\pi_2:C\times C_{k-2}\rightarrow C_{k-2}$. We obtain that
$$c_1(\F_j)=\Bigl[(\pi_2)_*\Bigl(\bigl(1+j((3k-6)\eta-\gamma-x)+\frac{j^2}{2}((3k-6)\eta-
\gamma-x)^2\bigr)\cdot\bigl(1-(2k-4)\eta\bigr)\Bigr)\Bigr]_1=$$
$$=-j^2\theta-j(3j-2)(k-2)x.$$
\end{proof}

Next we compute the class of $X_{j, \alpha}$ when $2\leq j\leq
k-1$. It is convenient to state our result as follows:

\begin{proposition}\label{xj}

Let $C$ be a general curve of genus $2k-2-j$ such that $2\leq
j\leq k-1$. We have the following relation in
$H^{4k-4-2j}\bigl(C\times \rm{Pic}$$^{3k-j-4}(C)\bigr)$:

$$\sum_{\alpha=0}^{[j/2]} a(j, \alpha) \ X_{j, \alpha}\equiv
\frac{2(2k-3)!}{k!(k-2)!(2k-3-j)!}\Bigl(\frac{\theta^{2k-2-j}}{2k-2-j}+\frac{(3kj+12k-9j-18)\theta^{2k-3-j}
\gamma}{2k-3} \ + $$
$$\frac{\theta^{2k-3-j}
\eta}{(2k-5)(2k-3)}
\bigl(60-64k+16k^2+j(60-\frac{67}{2}k-10k^2+6k^3)+j^2(15+2k-5k^2)+\frac{3}{2}
kj^3\bigr)\Bigr).
$$
\end{proposition}
\begin{proof} We fix an integer $0\leq \alpha\leq [j/2]$ and a divisor $D\in C_e$ of degree $e\geq k-j-1$.
We set $d:=3k-4-j$, $\Gamma:=D\times \mbox{Pic}^d(C)$  and denote
by $u, v:C\times C\times \mbox{Pic}^d(C)\rightarrow C\times
\mbox{Pic}^d(C)$ the two projections and by $J_{\alpha}(\L)$ the
$\alpha$-th jet bundle of the Poincare bundle on $C\times
\mbox{Pic}^d(C)$. Then $\mathcal{P}(D):=(\pi_2)_*\bigl(\L\otimes
\OO_{C\times \mbox{Pic}^d(C)}(\Gamma)\bigr)$ is a vector bundle of
rank $e+k-1$. For each $\alpha\geq 0$ we define the vector bundle
$J_{\alpha}(\L, D):=u_*\bigl(v^*(\L)\otimes \OO_{(\alpha+1)
\Delta+v^*(\Gamma)}\bigr)$,  where $\Delta\subset C\times C\times
\mbox{Pic}^d(C)$ is the diagonal, and consider the vector bundle
map
$$\phi_{\alpha}:J_{\alpha}(\L,
D)^{\vee}\rightarrow (\pi_2)^*(\mathcal{P}(D))^{\vee},$$ which on
fibres is the dual of the map $H^0(L_C\otimes \OO(D))\rightarrow
H^0(L_C\otimes \OO_{(\alpha+1)y+D}(D))$ for each pair $(y, L_C)\in
C\times \mbox{Pic}^d(C)$.

The cycle $X_{j, \alpha}$ consists of pairs $(y, L_C)$ such that
$\mbox{dim}\mbox{ Ker}(\phi_{\alpha}(y, L_C))\geq 1$ and
$\mbox{dim} \mbox{ Ker}(\phi_{j-\alpha+1} (y, L_C))\geq 2$ and its
class is given by the formula (cf. \cite{Fu},Theorem 14.3):
$$X_{j, \alpha}\equiv c_{j+1-\alpha}(\alpha) \
c_{k-1+\alpha-j}(j-\alpha+1)-c_{k-\alpha}(\alpha)\
c_{k-2+\alpha-j}(j-\alpha+1),$$ where
$c_l(\beta):=c_l\bigl((\pi_2)^*(\mathcal{P}(D))^{\vee}-J_{\beta}(\L,
D)^{\vee}\bigr)$. Since
$c_{tot}\bigl((\pi_2)^*(\mathcal{P}(D)\bigr)=e^{\theta}$ and
$c_1(\L)=d \eta+\gamma$ (independent of $D$!), after a short
calculation we obtain that
$$c_l(\beta)=\frac{\theta^l}{l!}+\frac{\theta^{l-1}}{(l-1)!}\bigl((\beta+1)\gamma+(\beta+1)(d+\beta(2k-3-j))\eta\bigr)
-\frac{\theta ^{l-1}}{(l-2)!}(\beta+1)(\beta+2)\eta.$$ A
straightforward computation now leads to the stated formula.
\end{proof}

For each integers $0\leq a\leq k-2$ and $ b\geq 2$ we shall define
vector bundles $\G_{a,b}$ over $\widetilde{\mathfrak
G}^{k-2}_{3k-6}$ with fibre
$$\G_{a,b}(C,L)=H^0(\wedge^a M_L\otimes L^{\otimes b})$$ over every point
$(C, L)\in \mathfrak G^{k-2}_{3k-6}$. Note that $\G_{i, 2 |
\mathfrak G^{k-2}_{3k-6}}=\mathcal{B}$, where $\mathcal{B}$ is the
vector bundle we introduced in Proposition \ref{ni}. Of course,
the question is how to extend this description over the divisors
$\Delta_j^0$. First we will extend $\G_{a, b}$ over
$\sigma^{-1}(\cM_g^0\cup \Delta_0^0\cup \Delta_1^0)$ and we start
by constructing the vector bundles $\G_{0, b}$:

\begin{proposition} For each $b\geq 2$ there
exists a vector bundle $\G_{0,b}$ over $\sigma^ {-1}(\cM_g^0\cup
\Delta_1^0 \cup \Delta_0^0) \subset \widetilde{\mathfrak
G}^{k-2}_{3k-6}$ of rank $k(3b-2)-6b+3$ whose fibres admit the
following description:
\begin{itemize}
\item For $(C, L)\in \mathfrak G^{k-2}_{3k-6}$, we have that
$\G_{0,b}\bigl(C,L)=H^0(C, L^{\otimes b})$. \item For $t=(C\cup_y
E, L)\in \sigma^{-1}(\Delta_1^0)$, where $L$ is the linear series
$\mathfrak g^{k-2}_{3k-6}$ on $C$ determining a limit $\mathfrak
g^{k-2}_{3k-6}$ on $C\cup_y E$, we have that
$$\G_{0,b}(t)=H^0(C, L^{\otimes b}(-2y))+\mathbb C\cdot
u^b\subset H^0\bigr(C, L^{\otimes b}),$$ where $u\in H^0(C, L)$ is
any section such that $\rm{ord}$$_y(u)=0$. \item For $t=(C/y\sim
q, L)\in \sigma^{-1}(\Delta_0^0)$, where $q,y\in C$ and $L$ is a
$\mathfrak g^{k-2}_{3k-6}$ on $C$, we have that
$$\G_{0,b}(t)=H^0(C, L^{\otimes b}(-y-q))\oplus \mathbb C\cdot u^b\subset H^0(C, L^{\otimes b}),$$
where $u\in H^0(C, L)$ is a section such that $\rm{ord}$$_y(u)=\rm{ord}$$_q(u)=0$.
\end{itemize}
\end{proposition}
\begin{proof} Recalling that $\sigma: \widetilde{\mathfrak G}^r_d\rightarrow \widetilde{\cM}_g$
denotes the variety of limit $\mathfrak g^r_d$'s over all treelike
curves of genus $g$, the morphism $\mathfrak
G^{k-2}_{3k-6}\rightarrow \mathfrak G^{k(3b-2)-6b+2}_{b(3k-6)}$
given by $(C, L)\mapsto (C, L^{\otimes b})$ can be extended to a
morphism $\nu_b:\widetilde{\mathfrak G}_{3k-6}^{k-2}\rightarrow
\widetilde{\mathfrak G}^{k(3b-2)-6b+2}_{b(3k-6)}$ over $\cM_g\cup
\Delta_0^0\cup \Delta_1^0$ as follows. We fix a point
$t=\bigl(C\cup_y E, l=(l_C, l_E)\bigr)$ a limit $\mathfrak
g^{k-2}_{3k-6}$, where $l_C=(L_C, H^0(L_C))\in G^{k-2}_{3k-6}(C)$
is such that $H^0(L_C)=H^0(L_C(-2y))\oplus \mathbb C\cdot u$, for
a certain $u\in H^0(L_C)$ with $\mbox{ord}_y(u)=0$, while
$l_E=\bigl(\OO_E((3k-6)y), (2k-6)p+|V_E|\bigr)$ (any limit
$\mathfrak g^{k-2}_{3k-6}$ on $C\cup_y E$ is of this form). Then
we set $\nu_b(t):=(l_C^b, l_E^b)$, where $l_C^b=\bigl(L_C^{\otimes
 b}, H^0(L_C^{\otimes b}(-2y))\oplus \mathbb C\cdot u^b\subset
H^0(L_C^{\otimes b})\bigr)$ and $l_E^b=\bigl(\OO_E(b(3k-6)p),
(2k-3)p+|\bigl(k(3b-2)-6b+3\bigr)p|\bigr)$. In other words,
$(l_C^b, l_E^b)$ is the $b$-th power of the original limit linear
series $\mathfrak g_{3k-6}^{k-2}$.

\noindent For a point $z=(C/y\sim q, L_C)\in \Delta_0^0$ with
$L_C$ being a $\mathfrak g^{k-2}_{3k-6}$ such that
$H^0(L_C)=H^0(L_C(-y-q))\oplus \mathbb C\cdot u$, we define
$\nu_b(z):=\bigl(L_C^{\otimes b}, H^0(L_C^{\otimes b}(-y-q))\oplus
\mathbb C\cdot u^b\subset H^0(L_C^{\otimes b})\bigr)$. The fact
that $\nu_b$ can be constructed algebraically follows easily from
the equations of the scheme $\widetilde{\mathfrak G}^{k-2}_{3k-6}$
described in \cite{EH1}, pg. 358. The variety
$\widetilde{\mathfrak G}_{b(3k-6)}^{k(3b-2)-6b+2}$ carries a
tautological vector bundle $\mathcal{T}=\mathcal{T}_C$ with fibre
over each point corresponding to a curve from $\cM_g^0\cup
\Delta_0^0$ being the space of global sections of the linear
series, while the fibre over a point corresponding to a curve from
$\Delta_1^0$ is the space of sections of the aspect of the limit
linear series corresponding to the curve of genus $g-1$. We define
$\G_{0, b}:=\nu_b^*(\mathcal{T})$. The description of the fibres
of the vector bundle $\G_{0, b}$ is now immediate.
\end{proof}

Having defined the vector bundles $\G_{0, b}$ we now define
inductively all vector bundles $\G_{a, b}$: first we define
$\G_{1, b}$ as the kernel of the multiplication map $\G_{0,
1}\otimes \G_{0, b}\rightarrow \G_{0, b+1}$, that is, by the exact
sequence
$$0\longrightarrow \G_{1, b}\longrightarrow \G_{0,1}\otimes \G_{0, b}\longrightarrow \G_{0, b+1}\longrightarrow 0.$$
Having defined $\G_{l, b}$ for all $l\leq a-1$, we define the
vector bundle $\G_{a, b}$ by the exact sequence
\begin{equation}\label{gi}
0\longrightarrow \G_{a, b}\longrightarrow \wedge^a \G_{0,
1}\otimes \G_{0, b}\stackrel{d_{a, b}}\longrightarrow \G_{a-1,
b+1}\longrightarrow 0.
\end{equation}
We now prove the right-exactness of the sequence (\ref{gi}) which will ensure the correctness of the definition of $\G_{a, b}$:  
\begin{proposition}
For  integers $b\geq 2$ and $1\leq a\leq i$, the Koszul
multiplication map $d_{a,b}:\wedge^a \G_{0, 1}\otimes \G_{0,
b}\rightarrow \G_{a-1, b+1}$ is well-defined and surjective. In
particular one can define the vector bundle $\G_{a,
b}:=\rm{Ker}$$(d_{a, b})$ over $\sigma^{-1}(\cM_g^0\cup
\Delta_0^0\cup \Delta_1^0)$ and the sequence $(\ref{gi})$ makes
sense.
\end{proposition}
\begin{proof} We do induction on $a$ and start with the case $a=1$.
 We only check the surjectivity
of the map $\G_{0, 1}\otimes \G_{0, b}\rightarrow \G_{0, b+1}$
over the locus $\sigma^{-1}(\Delta_1^0)$, the other cases, namely
$\sigma^{-1}(\cM_g^0)$ and $\sigma^{-1}(\Delta_0^0)$ being quite
similar. It is enough to show that if $C$ is a sufficiently
general curve of genus $2k-3$, $y\in C$ is a point and $L\in
W^{k-2}_{3k-6}(C)$ is a linear system with a cusp at $y$, then the
map $H^0(L)\otimes H^0(L^{\otimes b}(-2y))\rightarrow
H^0(L^{\otimes (b+1)}(-2y))$ is onto. This is equivalent to
$H^1(M_L\otimes L^{\otimes b}(-2y))=0$ which follows because of Proposition \ref{vanishing}.

We now treat the general case  $1\leq a\leq i$ and we want to show
that the Koszul map $d_{a, b}:\wedge^a H^0(L)\otimes \G_{0,
b}(L)\rightarrow \G_{a-1, b+1}(L)$ is surjective for each $b\geq
2$ and for each $(C,L)$ as above. For simplicity we have denoted
$\G_{a, b}(L)=\G_{a, b}(t)$, where $t\in \sigma^{-1}(\Delta_1^0)$
is the point $(C\cup_y E, l)$, $E$ is an arbitrary  elliptic tail
and $l$ is a limit linear series having $L$ as its $C$-aspect. We
can then write $H^0(L)/H^0(L(-2y))=\mathbb C\cdot u$, where $u\in
H^0(L)$ is uniquely determined up to translation by an element
from $H^0(L(-2y))$. We first claim that the restriction of the
Koszul differential
$$d_{a, b}^0:\wedge^a H^0(L)\otimes H^0(L^{\otimes
b}(-2y))\rightarrow H^0\bigl(\wedge^{a-1} M_L\otimes L^{\otimes
(b+1)}(-2y)\bigr)$$ is surjective. This is so because
$\mbox{Coker}(d_{a,b}^0)=H^1(\wedge^{a}M_L\otimes L^{\otimes
b}(-2y))=0$ (use again that $M_L$ is a stable bundle). But then
since $$H^0\bigl(\wedge^{a-1} M_L\otimes L^{\otimes
(b+1)}(-2y)\bigr)\subset \G_{a-1, b+1}(L)$$ is a linear subspace
of codimension ${k-2\choose a-1}$, to prove that $d_{a, b}$ itself
is surjective it suffices to notice that the image of the
injective map $\wedge^{a-1} H^0(L(-2y))\rightarrow \wedge^{a-1}
H^0(L)\otimes \G_{0, b+1}(L)$ given by
$$ f_1\wedge \ldots \wedge f_{a-1}\mapsto
f_1\wedge \ldots \wedge
 f_{a-1}\otimes u^{b+1}-\sum_{l=1}^{a-1} (-1)^{l-1} u\wedge
 f_1\wedge \ldots
\wedge \check{f_l} \wedge \ldots \wedge f_{a-1}\otimes f_l u^b$$
is entirely contained in $\mbox{Im}(d_{a, b})$ and is clearly
disjoint from $H^0\bigl(\wedge^{a-1} M_L\otimes L^{\otimes
(b+1)}(-2y)\bigr)$.
\end{proof}

To compute the intersection numbers of the divisors
$\overline{\mathcal{Z}}_{g, i}$ with the test curves $C^0$ and $C^1$
we need to understand the restriction of the vector bundle $\G_{i,
2}$ to $X$ and $Y$.

\begin{proposition}\label{x}
Let $C$ be a general curve of genus $2k-3$ and $k=3i+6$ with $i\geq
0$. If $C^1$ is the test curve in $\Delta_1$ obtained by attaching
to $C$ a fixed elliptic tail at a varying point of $C$ and
$X=\sigma^{-1}(C^1)$, then we have the following formula:
$$c_1(\G_{{2, i} |X})={3i+4\choose i}(c_x\ x+c_{\eta}\ \eta+c_{\gamma}\ \gamma+c_{\theta}\ \theta),$$
where
$$c_{\eta}=-27i+40, \mbox{ }c_{x}=-\frac{27i^4+153i^3+331i^2+323i+120}{(3i+4)(i+1)}$$
$$c_{\gamma}=\frac{5i+8}{3i+4}, \mbox{ }c_{\theta}=-\frac{27i^3+101i^2+124i+48}{(3i+3)(3i+4)}.$$
\end{proposition}

\begin{proof} For $j\geq 2$ we define the jet bundle
$J_1(\cM^{\otimes j}):=(u_*)\bigl(v^*(\OO_{C\times
C_{k-2}}(\cM^{\otimes j}))\otimes \OO_{2 \Delta}\bigr)$, where $u,
v:C\times C\times C_{k-2}\rightarrow C\times C_{k-2}$ are the
projections. We have two exact sequences of vector bundles on $X$:
$$0\longrightarrow u_*\bigl(v^*(\cM^{\otimes j})\otimes
\OO(-2\Delta)\bigr)_{|X}\longrightarrow \G_{0, j
|X}\longrightarrow \cM^{\otimes j} _{|X}\longrightarrow 0$$ and
$$0\longrightarrow u_*\bigl(v^*(\cM^{\otimes j})\otimes
\OO(-2\Delta)\bigr)_{|X}\longrightarrow
\pi_2^*(\F_j)_{|X}\longrightarrow J_1(\cM^{\otimes
j})_{|X}\longrightarrow 0,$$ from which we can write that
$$c_1(\G_{0, j\ |X})=c_1(\pi_2^*(\F_j)_{| X})-c_1\bigl((\pi_1^*K_C\otimes \cM^{\otimes j})_{|X}\bigr)=$$
$$=-j^2\theta-(3j+4)(k-2)\ \eta+j\gamma-j((3j-2)(k-2)-1)\ x.$$
From the exact sequences defining $\G_{i, 2}$ we then obtain that
$$c_1(\G_{i, 2 |X})=\sum_{l=0}^i (-1)^l c_1(\wedge^{i-l}\G_{0,
1 |X}\otimes \G_{0, l+2 |X}).$$ Note that Riemann-Roch gives that
$\mbox{rk}(\G_{0, l+2})=(3k-6)l+4k-9$,  while
$c_1(\wedge^{i-l}\G_{0, 1 |X})={k-2\choose
i-l-1}(-\theta-(k-2)x).$ To obtain a closed formula for
$c_1(\G_{i, 2 |X})$ we now specialize to the case $k=3i+6$ and we
write that
$$c_1(\G_{i, 2 |X})=\sum_{l=0}^i
(-1)^l\Bigl[\bigl((9i+12)l+(12i+15)\bigr){3i+4\choose
i-l-1}(-\theta-(3i+4)x)+$$ $$+ {3i+5\choose
i-l}\Bigl(-(l+2)^2\theta+(l+2)\gamma-(3i+4)(3l+10)\eta-(l+2)((3l+4)(3i+4)-1)x\Bigr)\Bigr].
$$
A long but elementary calculation leads to the stated formula.

\end{proof}

A similar calculation yields the first Chern class of $\G_{{i,
2} |Y}$:
\begin{proposition}\label{y}
Let $(C, q)$ be a general pointed curve of genus $2k-3$ and $C^0$
the test curve in $\Delta_0$ obtained by identifying the fixed
point $q$ with a varying point $y$ on $C$. If
$Y=\sigma^{-1}(C^0)$, then we have the formula
$$c_1(\G_{{i, 2} |Y})={3i+4\choose i}(d_x\ x+ \eta+d_{\theta}\ \theta),$$
where
$$d_{x}=-\frac{27i^4+153i^3+331i^2+323i+120}{(3i+4)(i+1)},\ d_{\theta}=-\frac{27i^3+101i^2+124i+48}{(3i+3)(3i+4)}. $$
\end{proposition}
\begin{proof} We define $\Gamma_q:=C\times \{q\}\times C_{k-2}$ and
denote by $\iota:C\times C_{k-2}\hookrightarrow \Gamma_q$ the
inclusion map. Then we have the exact sequences on $Y$
$$0\longrightarrow (u_*)\bigl(v^*(\cM^{\otimes
j})\otimes \OO(-\Delta-\Gamma_q)\bigr)_{|Y}\longrightarrow \G_{0,
j|Y}\longrightarrow \cM^{\otimes j}_{|Y}\longrightarrow 0,$$ and
$$0\longrightarrow \iota^*\bigl(v^*(\cM^{\otimes j})\otimes
\OO(-\Delta)\bigr)_{|Y} \longrightarrow
(u_*)\bigl(v^*(\cM^{\otimes j})\otimes
\OO_{\Delta+\Gamma_q}\bigr)_{|Y}\longrightarrow
\pi_2^*(\F_j)_{|Y}\longrightarrow 0.$$ Using that $c_1(\cM)=
(3k-6)\eta-\gamma-x$, it is now straightforward to check that
$$c_1\bigl(\iota^*(v^*(\cM^{\otimes j})\otimes
\OO(-\Delta))\bigr)=-jx-\eta,$$ hence $c_1(\G_{{0, j}
|Y})=c_1(\pi_2^*(\F_j))_{|Y}+jx+\eta=-j^2\theta+\eta-j((3j-2)(k-2)-1)x.$
But then as in the previous proposition we have that
$$c_1(\G_{{i, 2} |Y})=\sum_{l=0}^i (-1)^l c_1(\wedge^{i-l}\G_{{0,
1}|Y}\otimes \G_{{0, l+2} |Y}),$$ which after some calculations
yields the stated formula in the special case $k=3i+6$.
\end{proof}

In what follows we extend the vector bundles $\G_{a, b}$ over the
whole $\widetilde{\mathfrak G}_{3k-6}^{k-2}$. Roughly speaking, the
fibre of $\G_{a, b}$ over a point corresponding to a singular curve
$C\cup_y D$  will be the space $H^0(C\cup D, \wedge^a M_{L}\otimes
L^{\otimes b})$, where $L$ will be a line bundle on $C\cup D$
obtained by twisting appropriately a limit $\mathfrak
g^{k-2}_{3k-6}$ on $C\cup D$. The twisting is chosen such that, the
vector bundles $\G_{a, b}$ sit in the exact sequences (\ref{gi}). In
the next statement we use the notation introduced in Proposition
\ref{limitlin3}:

\begin{theorem}\label{dj}
For integers $0\leq a \leq i$, $ b\geq 2$, or $(a, b)=(0, 1)$,
there exists a vector bundle $\G_{a, b}$ over
$\widetilde{\mathfrak G}^{k-2}_{3k-6}$ having the following
properties:

\noindent (1) If $(C, L)\in \mathfrak G^{k-2}_{3k-6}$, then
$\G_{a, b}(C, L)=H^0(\wedge^a M_L\otimes L^{\otimes b})$, that is,
$\G_{a, b \mbox{ }| \mathfrak G^{k-2}_{3k-6}}=\G_{a, b}$.

If $t=(C\cup_y D, L_C, L_D) \in \sigma^{-1}(\Delta_j^0)\subset
\widetilde{\mathfrak G}^{k-2}_{3k-6}$, with $2\leq g(D)=j\leq k-1,
g(C)=2k-2-j, L_C\in W^{k-2}_{3k-4-j}(C)$ and $L_D\in
W^{k-2}_{k-2+j}(D)$, then there are two situations:

\noindent (2) If $2\leq j\leq 2i+2$, then $\G_{a, b}(t)=H^0(C\cup
D, \wedge^a M_{L_{C\cup D}}\otimes L_{C\cup D}^{\otimes b}),$
where $L_{C\cup D}=\bigl(L_C(-(j+2)y), L_D(-(k-2-j)y)\bigr)\in
\rm{Pic}$$^{3k-6-2j}(C)\times \rm{Pic}$$^{2j}(D)$ is a globally
generated line bundle on $C\cup D$.

\noindent (3) If $2i+3\leq j\leq 3i+5$ and $c:=[(j+2)/2]$, then
$\G_{a, b}(t)=H^0(C\cup D, \wedge^a M_{L_{C\cup D}}\otimes
L_{C\cup D}^{\otimes b})$, where $L_{C\cup D}=\bigl(L_C(-cy)),
L_D(-(k-c)y)\bigr)\in \rm{Pic}$$^{3k-4-j-c}(C)\times
\rm{Pic}$$^{j+c-2}(D)$ is a globally generated line bundle on
$C\cup D$.

\noindent (4) There are exact sequences over $\widetilde{\mathfrak
G}^{k-2}_{3k-6}$:\mbox{ } $0\longrightarrow \G_{a,
b}\longrightarrow \wedge^a \G^{0, 1}\otimes \G_{0,
b}\longrightarrow \G_{a-1, b+1}\longrightarrow 0.$

\end{theorem}
\begin{proof} For $j\geq 2$ the divisors $\Delta_j ^0$ are mutually
disjoint and we can carry out the construction of $\G_{a, b}$ over
each open set $\sigma^{-1}(\cM_g^0\cup \Delta_j^0)\subset
\widetilde{\mathfrak G}^{k-2}_{3k-6}$ and glue the resulting
bundles together. We denote by $\pi:\mathcal{C}_j\rightarrow
\cM_g^0\cup \Delta_j^0$ the restriction of the universal curve and
by $p:\mathcal{Y}=\mathcal{C}_j \times _{\cM_g^0\cup \Delta_j^0}
\widetilde{\mathfrak G}_{3k-6}^{k-2}\rightarrow
\widetilde{\mathfrak G}^{k-2}_{3k-6}$. Then $(\sigma
p)^{-1}(\Delta_j^0)=D_{j}+D_{2k-2-j}$, where $D_j$ (resp.
$D_{2k-2-j}$) denotes the divisor in $\mathcal{Y}$ corresponding
to the marked point being on the genus $j$ (resp. $2k-2-j$)
component. If $\L=\L_{0, 3k-6}$ is the Poincare bundle on
$\mathcal{Y}$ chosen such that it parametrizes bundles having
bidegree $(0, 3k-6)$  on curves $[D\cup_y C] \in \Delta_j^0$, then
for $b\geq 1$ we set
$$\G_{0, b}:=p_*\bigl(\L^{\otimes b}\otimes
\OO_{\mathcal{Y}}(-2bj \ D_j)\bigr), \mbox{ when }2\leq j\leq
2i+2,$$ and $\G_{0, b}:=p_*\bigl(\L^{\otimes b}\otimes
\OO_{\mathcal Y}(-b(j-2+c) D_j)\bigr)$ when $2i+3\leq j\leq 3i+5$.
For each point $t=(R=C\cup_y D, L_C, L_D) \in
\sigma^{-1}(\Delta_j^0)$, we have that $\G_{0, b}(t)=H^0(R,
L_R^{\otimes b})$, where in the case $j\leq 2i+2$ we have
$L_R=(L_C(-(j+2)y), L_D(-(k-2-j)y))\in \mbox{Pic}^{3k-6}(R)$ and
since $H^1(R, L_R^{\otimes b})=0$ for $b\geq 2$, it follows that
$\G_{0, b}$ is a vector bundle for all $b\geq 1$. When $1\leq
a\leq i$, the sheaves $\G_{a, b}$ are defined inductively using
the sequences (\ref{gi}). For this to make sense and in order to
conclude that $\G_{a, b}$ are vector bundles, one has to show that
\begin{equation}\label{vanr}
H^1(R, \wedge^a M_{L_R}\otimes L_R^{\otimes b})=0,
\end{equation}
for all $b\geq 2$ and $0\leq a\leq i$. To achieve this we use the
Mayer-Vietoris sequence on $R$:
$$0\rightarrow \wedge^a M_{L_R}\otimes L_R^{\otimes
b}\longrightarrow \bigl(\wedge^a M_{L_R}\otimes L_R^{\otimes
b}\bigr)_{|C}\oplus \bigl(\wedge^a M_{L_R}\otimes L_R^{\otimes
b}\bigr)_{| D}\longrightarrow \wedge^a M_{L_R}\otimes L_R^{\otimes
b}\otimes \mathbb C(y)\rightarrow 0,$$ together with the exact
sequences
$$0\longrightarrow  H^0(L_{R|C}(-y))\otimes
\OO_D\longrightarrow M_{L_{R}}\otimes \OO_D\longrightarrow
M_{L_{R|D}}\longrightarrow 0,$$ and
$$0\longrightarrow H^0(L_{R|D}(-y))\otimes \OO_C\longrightarrow
M_{L_R}\otimes \OO_C\longrightarrow M_{L_{R|C}}\longrightarrow
0.$$ We obtain that (\ref{vanr}) holds if we can show that
$$H^1(C, \wedge^a M_{L_{R|C}}\otimes L_{R|C}^{\otimes b})=H^1(C,
\wedge^a M_{L_{R|C}}\otimes L_{R|C}^{\otimes b}(-y))=0$$ and $$H^1(D, \wedge^a M_{L_{R|D}}\otimes L_{R|D}^{\otimes
b})=H^1(D, \wedge^a M_{L_{R|D}}\otimes L_{R|D}^{\otimes
b}(-y))=0.$$ We only check this when $j\leq 2i+2$ the remaining
case being similar. Since $h^0(L_{R|C})=k-1-j$ and
$h^0(L_{R|D})=j+1$, the inequalities
$g(C)+\mbox{min}\{h^0(L_{R|C})-1, i\}\leq \mbox{deg}(L_{R|C})$ and
$g(D)+\mbox{min}\{h^0(L_{R|D})-1, i\}\leq \mbox{deg}(L_{R|D})$ are
satisfied precisely because $j\leq 2i+2$ and we can invoke Lemma
\ref{vanishing}.
\end{proof}

Next, for $0\leq a\leq i$ and $b\geq 1$ we define vector bundles
$\H_{a, b}$ over $\widetilde{\mathfrak G}_{3k-6}^{k-2}$ having
fibre $\H_{i, j}(C, L)=H^0(\wedge^i M_{\PP^{k-2}}(j))$ over each
point $(C,L)\in \mathfrak G^{k-2}_{3k-6}$ corresponding to an
embedding $C\stackrel{|L|}\hookrightarrow \PP^{k-2}$. First we set
$\H_{0, 1}:=\G_{0, 1}$ and $\H_{0, b}:=\mbox{Sym}^b \H_{0, 1}$ for $b\geq 1$. 
Having already defined $\H_{a-1, b}$ for all $b\geq 1$,  we define $\H_{a, b}$
via the exact sequence
\begin{equation}\label{sym}
0\longrightarrow \H_{a, b}\longrightarrow \wedge^a \H_{0,
1}\otimes \mbox{Sym}^b \H_{0, 1}\longrightarrow \H_{a-1,
b+1}\longrightarrow 0.
\end{equation}
The bundles $\H_{a, b}$ being defined entirely in terms of $\H_{0, 1}$, the right exactness of the sequence (\ref{sym}) is obvious.
There is a natural vector bundle morphism $\phi_{a, b}:\H_{a,
b}\rightarrow \G_{a, b}$ over $\widetilde{\mathfrak
G}^{k-2}_{3k-6}$; when $k=3i+6, a=i$ and $b=2$, then
$\mbox{rank}(\H_{i, b})=\mbox{rank}(\G_{i, b})$ and the
degeneracy locus of $\phi_{i, 2}$ is the codimension $1$
compactification over $\widetilde{\mathfrak G}_{3k-6}^{k-2}$ of
the locus $\mathcal{U}_{g, i}$ defined in Proposition \ref{ni}.

Next we compute the Chern classes of $\H_{a, b}$ along the curves
$X$ and $Y$:
\begin{proposition}\label{has}
When $k=3i+6$, we  have the following formulas for the first Chern
class of $\H_{i, 2}$:
$$c_1(\H_{i, 2 |X})=c_1(\H_{i, 2 |Y})=-(3i+6){3i+4\choose i}\bigl((3i+4)x+\theta\bigr).$$
\end{proposition}
\begin{proof} Using repeatedly the sequence (\ref{sym}) we obtain
the  formula
$$c_1(\H_{i, 2 |X})=c_1(\H_{i, 2 |Y})=\sum_{l=0}^i (-1)^l c_1(\wedge ^{i-l}\G_{0, 1
|X}\otimes \mbox{Sym}^{l+2} \G_{0, 1 |X}).$$
Since $c_1(\G_{0, 1})=-\theta-(3i+4)x$ and $c_1(\mbox{Sym}^j(\G_{0, 1}))=\frac{j(3i+4+j)}{(3i+5)(3i+4)}{3i+j+3\choose j}
 c_1(\G_{0, 1})$, the stated formula is obtained after a straightforward calculation.
\end{proof}

In the remainder of this section we will compute $c_1(\G_{i,
2}-\H_{i, 2})\cdot \sigma^*(C^j)$, where $C^j\subset \Delta_j$ is
the test curve associated to $[C]\in \cM_{2k-2-j}$. We retain the
notation introduced in Propositions \ref{limitlin2} and
\ref{limitlin3}. There are two definitions of the bundles $\G_{a,
b |\sigma^*(\Delta^0_j)}$ depending whether $j\leq 2i+2$ or $j\geq
2i+3$ and we will explain in detail the calculations only in the
first case, the second being similar. We start by describing the
Chern number of the restriction of $\G_{0, 1}$ to the $1$-cycle
$\sum_{\alpha=0}^{[j/2]} a(j, \alpha)\ X_{j, \alpha}$ on
$C\times \mbox{Pic}^{3k-j-4}(C)$:
\begin{proposition}\label{xja}
Suppose $C$ is a general curve of genus $2k-2-j$ where $2\leq
j\leq 2i+2$ and let $C^j\subset \Delta_j\subset \mm_{2k-2}$ be the test curve
associated to $C$. Then
$$c_1(\G_{0,1 | \sum_{\alpha=0}^{[j/2]} a(j, \alpha)
X_{j,
\alpha}})=-\theta-j\bigl((3k-j-4)\eta+\gamma\bigr)-j(j+1)(2k-3-j)
\eta+\aleph,$$
$$ \mbox{ where }\mbox{ } \mbox{ }\  \aleph= 6(j+1)\frac{{6i+10-j\choose 3i+6}{6i+8\choose j}}{{3i+4\choose j}}.$$
\end{proposition}
\begin{proof} For each $0\leq \alpha \leq j/2$, we fix once and for all a linear series $L_D\in
W^{k-2}_{k-2+j}(D)$ such that $h^0(L_D(-(k+1-\alpha)y))=\alpha$
and $h^0(L_D(-(k+\alpha-j)y))=j-\alpha$ (there are $a(j, \alpha)$
such $\mathfrak g^{k-2}_{k-2+j}$'s and for each of them $X_{j,
\alpha}\times \{L_D\}$ is a component of $\sigma^*(C^j)$). Just
like in Proposition \ref{xj}, we introduce the Poincare bundle
$\L$ on $X_{j, \alpha}$ and we define the sheaf on $X_{j, \alpha}$
given by $\mathcal{P}:=\pi_2^*(\pi_{2 * }\L)$. Furthermore, for each
$a\geq 0$ we consider the Taylor map $\nu_a:
\mathcal{P}\rightarrow J_a(\L)$ (the target being the $a$-th jet
bundle of $\L$), and we set $\F_a:=\mbox{Im}(\nu_a)\subset
J_a(\L)$ and $\mathcal{P}_{a+1}:=\mbox{Ker}(\nu_a)$.

For each pair $(y, L_C)\in X_{j, \alpha}$ we have an exact sequence
$$0\longrightarrow H^0(D, L_{R|D}(-y))\longrightarrow H^0(C\cup_y D,
L_R)\longrightarrow H^0(C, L_{R|C})\longrightarrow 0,$$ where
 $L_R=\bigl(L_C(-(j+2)y), L_D(-(k-2-j)y)\bigr)$ is the line
bundle on $C\cup D$ whose global sections give the fibre of
$\G_{0, 1}$ at each point in $\widetilde{\mathfrak
G}^{k-2}_{3k-6}$. This exact sequence globalizes to an exact
sequence of vector bundles on $X_{j, \alpha}$
$$0\longrightarrow H^0(L_D(-(k-1-j)y))\otimes \OO_{X_{j,
\alpha}}\longrightarrow \G_{0, 1 |X_{j, \alpha}}\longrightarrow
\mathcal{P}_{j+2}\longrightarrow 0,$$  thus $c_1(\G_{0,1|X_{j,
\alpha}})=c_1(\mathcal{P}_{j+2})$. We now study how the bundles
$\F_a$ relate to one another via the exact sequences
$0\longrightarrow \pi_1^*(K_C^{\otimes a})\otimes
\L\longrightarrow J_a(\L)\longrightarrow
J_{a-1}(\L)\longrightarrow 0$, \  linking successive jet bundles.
Recalling that for a generic point $(y, L_C)\in X_{j, \alpha}$
we have that $a^{L_C}(y)=(0, \ldots, \alpha-1, \alpha+1, \ldots ,
j-\alpha, j-\alpha+2, \ldots, k)$ (cf. Proposition
\ref{limitlin2}), it follows that $\F_a=J_a(\L)$ for $a\leq
\alpha-2$, while $\F_{\alpha-1}$ is obtained from
$J_{\alpha-1}(\L)$ via an elementary transformation along the
divisor $D_1\subset X_{j, \alpha}$ consisting of the $b(2k-2-j,
k-2+\alpha-j)$ points $(y, L_C)\in X_{j, \alpha}$ where
$a_{\alpha-1}^{L_C}(y)=\alpha$. More precisely, one has the exact
sequence
$$0\longrightarrow \pi_1^*(K_C^{\otimes (\alpha-1)})\otimes
\L\otimes \OO_{X_{j, \alpha}}(-D_1)\longrightarrow
\F_{\alpha-1}\longrightarrow J_{\alpha-2}(\L)\longrightarrow 0.$$
Furthermore, we have exact sequences
$$0\longrightarrow \pi_1^*(K_C^{\otimes (\alpha+1)})\otimes
\L\longrightarrow \F_{\alpha+1}\longrightarrow
\F_{\alpha-1}\longrightarrow 0\ \mbox{ and }$$
$$0\longrightarrow \pi_1^*(K_C^{\otimes a})\longrightarrow
\F_a\longrightarrow \F_{a-1}\longrightarrow 0, \  \mbox{ for all }
\alpha+2\leq a\leq j-\alpha-2.$$ If $D_2\subset X_{j, \alpha}$ is
the divisor consisting of the $b(2k-2-j, k-1+\alpha-j)$ points
$(y, L_C)$ satisfying the condition
$a^{L_C}_{j-\alpha-1}(y)=j-\alpha+1$, then we also have the exact
sequences
$$0\longrightarrow \pi_1^*(K_C^{\otimes (j-\alpha)})\otimes
\L\otimes \OO_{X_{j, \alpha}}(-D_2)\longrightarrow
\F_{j-\alpha}\longrightarrow \F_{j-\alpha-1}\longrightarrow 0,$$
$$0\longrightarrow \pi_1^*(K_C^{\otimes (j-\alpha+2)})\otimes
\L\longrightarrow \F_{j-\alpha+2}\longrightarrow
\F_{j-\alpha}\longrightarrow 0,$$ and finally
$$0\longrightarrow \pi_1^*(K_C^{\otimes a})\otimes
\L\longrightarrow \F_{a}\longrightarrow \F_{a-1}\longrightarrow 0,
\  \mbox{ for  } j-\alpha+2\leq a\leq k-1.$$

\noindent We gather all the information contained in these
sequences to obtain that

$c_1(\F_{j+1})=j c_1(\L)+j(j+1)(g(C)-1)\eta-\mbox{deg}(D_1+D_2),$
\mbox{  } and then
$$c_1(\G_{0,1 |\sum_{\alpha} a(j, \alpha) X_{j,
\alpha}})=c_1(\mathcal{P})-c_1(\F_{j+1})=-\theta-j((3k-j-4)\eta+\gamma)-j(j+1)(2k-3-j)\eta+\aleph,$$
 $$\mbox{ where }\aleph:=\sum_{\alpha=0}^{[j/2]} a(j, \alpha)\bigl(b(2k-2-j,
k-2+\alpha-j)+b(2k-2-j, k-1+\alpha-j)\bigr)=$$
$$=j!(6i+10-j)!\sum_{\alpha=0}^{j+1}
\frac{(j-2\alpha-1)(j-2\alpha)(j-2\alpha+1)^2}{(3i+5-\alpha)!\
(3i+5+\alpha-j)!\ (j-\alpha+1)!\ \alpha!}.$$ This sum can be computed with
Maple which finishes the proof.
\end{proof}

 Next we compute the Chern class of $\G_{0, b}$
with $b\geq 2$. The proof being similar to that
 in Proposition \ref{xja} (in fact simpler), we decided to omit it.
\begin{proposition}\label{xjb}
Let $b\ge 2$ and $2\leq j\leq 2i+2$. Then for each $0\leq \alpha
\leq j/2$ we have that
$$c_1(\G_{0, b | X_{j, \alpha}})=-b^2\theta-b^2(j+2)(\gamma+(3k-4-j)\eta)-{bj+2b\choose 2}(2k-3-j)\eta.$$
\end{proposition}

We turn our attention to the remaining components of
$\sigma^*(C^j)$. First we show that $Y_{j, \alpha}^{'''}$ does not
appear in the computation of $c_1(\G_{i, 2}-\H_{i, 2})\cdot
\sigma^*(C^j)$:
\begin{proposition}\label{tri}
For $0\leq \alpha\leq j/2,  a\geq 0$ and
$b\geq 1$, we have that $c_1(\G_{a, b})\cdot Y_{j,
\alpha}^{'''}=c_1(\H_{a, b})\cdot Y_{j,
\alpha}^{'''}=0.$
\end{proposition}
\begin{proof} Clearly, it suffices to show that $c_1(\G_{0,
b})\cdot Y_{j, \alpha}^{'''}=0$ for all $b\geq 1$ and then use the
exact sequences (\ref{gi}) and (\ref{sym}). We carry this out only
for $b=1$ the case $b\geq 2$ being analogous. Fix $l=(L_D,
V_D\subset H^0(L_D))\in Y_{j, \alpha}^{'''}$. By duality
$K_D\otimes L_D^{\vee}((k+2-\alpha)y)$ is one of the $a(j,
\alpha)$ (that is, finitely many) linear systems $\mathfrak
g^1_{j-\alpha+1}$ on $D$ with a $(j-2 \alpha+1)$-fold point at
$y$. The space of sections $V_D\subset H^0(L_D)$ is chosen such
that $H^0(L_D(-2y))\subset V_D$, hence $Y_{j, \alpha}^{'''}$ can
be identified with the disjoint union of $a(j, \alpha)$ copies of
the projective line $\PP\bigl(H^0(L_D)/H^0(L_D(-2y))\bigr)$, one
for each choice of $L_D$. Then $H^0(L_{C\cup
D|D})=H^0(L_D(-(k-1-j)y))$, for each $l\in Y_{j, \alpha}^{'''}$
(that is, $\G_{0, 1}(l)$ is independent of $V_D$). Globalizing
this we get that $\G_{0, 1 | Y_{j, \alpha}^{'''}}$ is trivial.
\end{proof}

In a similar fashion we have the following:
\begin{proposition}\label{primul}
For $0\leq \alpha\leq j/2, a\geq 0$ and $b\geq 1$, we have the equalities $c_1(\G_{a, b})\cdot Y_{j, \alpha}^{'}=c_1(\H_{a, b})\cdot Y_{j, \alpha}^{'}=0.$
\end{proposition}

\begin{proposition}\label{secondy}
For $0\leq \alpha\leq j/2$ and $b\geq 1$, we have the equality
$$c_1(\G_{0,
b})\cdot Y_{j, \alpha}^{''}=-\frac{b^2 (j-2\alpha) j!}{\alpha!
(j-\alpha)!}.$$
\end{proposition}
\begin{proof}  We start by recalling that $Y_{j, \alpha}^{''}$ is
the locus of line bundles $L_D \in \mbox{Pic}^{k-2+j}(D)$ such that
$h^0(L_D(-(k+1+\alpha-j)y))\geq j-\alpha-1$ and
$h^0(L_D(-(k-\alpha+1)y))\geq \alpha$. Using \cite{Fu}, Theorem
14.3, one finds that  $Y_{j, \alpha}^{''}\equiv
\frac{j-2\alpha}{\alpha!(j-\alpha)!}\theta^{j-1}$. Having fixed once
and for all one of the finitely many linear series $L_C\in
W^{k-2}_{3k-j-4}(C)$ such that $\{L_C\}\times Y_{j,
\alpha}^{'}\subset \widetilde{\mathfrak G}^{k-2}_{3k-6}$, there is
an  identification $$\G_{0, 1}(L_D)=H^0(L_D(-(k-2-j))y),$$ for each
$L_D\in Y_{j, \alpha}^{''}$. Since the vector bundle
$\mbox{Pic}^{k-2+j}(D)\ni L_D\mapsto H^0(L_D\otimes \OO_{\beta}y)$,
is algebraically trivial for each $\beta \geq 0$ (remember that
$y\in D$ is a fixed point), we obtain that $c_1(\G_{0, 1 |  Y_{j,
\alpha}^{''}})=-\theta$. Similarly, for $b\geq 2$ we get that
$c_1(\G_{0, b | Y_{j, \alpha}^{''}})=-b^2\theta$ from which now the
conclusion follows.
\end{proof}

Propositions \ref{xja}, \ref{xjb}, \ref{tri}, \ref{primul},
\ref{secondy} will enable us to compute the intersection numbers
$c_1(\G_{i, 2}-\H_{i, 2})\cdot \sigma^*(C^j)$ needed  to determine
the coefficient of $\delta_j$ in the expansion of
$[\overline{\mathcal{Z}}_{6i+10, i}]$. Repeatedly using (\ref{sym})
, one can write that
\begin{equation}\label{hidoi}
c_1(\H_{i, 2})=\sum_{l=0}^i (-1)^l c_1\bigl(\wedge^{i-l} \G_{0,
1}\otimes \mbox{Sym}^{l+2} \G_{0, 1}\bigr)=(3i+6){3i+4\choose i}
c_1(\G_{0, 1})
\end{equation}
(note the similarity with Proposition \ref{has}). Using (\ref{gi})
we also have that
\begin{equation}\label{gidoi}
c_1(\G_{i, 2})=\sum_{l=0}^i (-1)^l c_1\bigl(\wedge^{i-l} \G_{0,
1}\otimes \G_{0, l+2}\bigr),\end{equation} and a simple calculation yields that
 \begin{equation}\label{diferenta}
c_1(\G_{i, 2}-\H_{i,
2})=-\frac{13i^2+35i+24}{(3i+4)(i+1)}{3i+4\choose i} c_1(\G_{0,
1})+\sum_{l=0}^i (-1)^l {3i+5\choose i-l}c_1(\G_{0, l+2}).
\end{equation}

\section{The slope of $\overline{\mathcal{Z}}_{6i+10, i}$}

In this section we finish the calculation of the class  of the virtual divisor $\zz_{g,
i}$. As before, we fix $i\geq 0, k=3i+6$ and $g=6i+10$, hence $\phi:\H_{i,
2}\rightarrow
\G_{i, 2}$ is a morphism between two vector bundles
of the same rank defined over $\tilde{\mathfrak G}^{3i+4}_{9i+12}$.

\begin{theorem}\label{formula}

We fix $i\geq 0$ and denote by $\sigma:\widetilde{\mathfrak G}^{3i+4}_{9i+12}\rightarrow \widetilde{\cM}_{6i+10}$ the
natural projection map. Then we have the following formula for the  class of the pushforward to $\widetilde{\cM}_{6i+10}$ of the virtual degeneracy locus of the morphism $\phi:\H_{i, 2}\rightarrow \G_{i, 2}$:
$$\sigma_*(c_1(\G_{i, 2}-\H_{i,, 2}))= \frac{1}{3i+5}{6i+7 \choose 3i+4, 2i+1, i+2}\Bigl(a \lambda-b_0\ \delta_0-\cdots
-b_{3i+5}\ \delta_{3i+5}\Bigr),$$ where
$$a=\frac{(4i+7)(6i^2+19i+12)}{(2i+3)(i+2)}, \mbox{ }\mbox{ }
b_0=\frac{12i^2+31i+18}{3(2i+3)},\mbox{ }  b_1=\frac{12i^2+33i+20}{i+2} \mbox{ and } $$
$$b_j=\frac{j f(i, j)}{6(i+1)(i+2)(2i+3)(6i+9-j)}, \mbox{ for } 2\leq j\leq 2i+2,$$
with
$$f(i, j)=864\ i^5-(240j-5256)i^4+(16j^2-998j+11830)i^3+$$
$$+(31j^2-1254j+11585)i^2
+(-5j^2-286j+3981)i-24j^2+240j-216,
$$
and
$$b_j=\frac{2(13i^2+35i+24)(3i+5)}{(i+1)(i+2){6i+10\choose j}}{3i+4\choose (j-1)/2}^2
+\frac{g(i, j)}{12(i+1)(i+2)(2i+3)(6i+9-j)}$$
 for an odd $j$ such that $2i+3\leq j\leq
3i+5$, where
$$g(i, j)=1728ji^5 -(576 j^2- 10512j+540)i^4 +(48j^3-2592j^2+24120j-3234)i^3+$$
$$(140j^3-3926j^2+25176j-7278)i^2+(109j^3-2107j^2+10875j-7263)i
+12j^3-156j^2+972j-2700$$
and
$$b_j=\frac{(13i^2+35i+24)((3i+5)(2j+1)-j^2)}{2(i+1)(i+2)(3i+5){6i+10\choose
j}}{3i+5\choose j/2}^2+\frac{h(i, j)}{12(i+1)(i+2)(2i+3)(6i+9-j)},$$
for an even $j$ such that $2i+3\leq j\leq 3i+5$, where
$$h(i, j)=1728ji^5 -(576 j^2- 10512j+1020)i^4 +(48j^3-2592j^2+24280j-6034)i^3+$$
$$(140j^3-3942j^2+25830j-13400)i^2+(109j^3-2145j^2+11774j-13230)i
+12j^3-180j^2+1392j-4896.$$

 In particular, $b_j\geq b_0$ for all $1\leq j\leq 3i+5$ and
 $$s(\sigma_*(c_1(\G_{i, 2}-\H_{i, 2})))=\frac{a}{b_0}=\frac{3(4i+7)(6i^2+19i+12)}{(12i^2+31i+18)(i+2)}.$$

\end{theorem}

\begin{proof}
Since $\mbox{codim}(\mm_g-\widetilde{\cM}_g, \mm_g)\geq 2$, it
makes no difference whether we compute the class $\sigma_*(\G_{i,
2}-\H_{i, 2})$ on $\widetilde{\cM}_g$ or on $\mm_g$ and we can
write
$$\sigma_*(\G_{i, 2}-\H_{i, 2})= A \lambda-B_0\ \delta_0-B_1\ \delta_1-\cdots
-B_{3i+5}\ \delta_{3i+5}\in \mbox{Pic}(\mm_g),$$
 where $\lambda, \delta_0, \ldots, \delta_{3i+5}$ are the
generators of $\mbox{Pic}(\mm_g)$. We start with the following:

\noindent {\bf Claim:} \emph{One has the relation
$A-12B_0+B_1=0$.}

We pick a general curve $[C, q]\in \cM_{6i+9, 1}$ and at the fixed
point $q$ we attach to $C$ a Lefschetz pencil of plane cubics. If
we denote by $R\subset \mm_g$ the resulting curve, then
 $R\cdot \lambda=1, \ R\cdot \delta_0=12, \ R\cdot
\delta_1=-1$ and $R\cdot \delta_j=0$ for $j\geq 2$. The
relation $A-12B_0+B_1=0$ follows once we show that
$\sigma^*(R)\cdot c_1(\G_{i, 2}-\H_{i, 2})=0$. To achieve this we
check that $\G_{0, b |\sigma^*(R)}$ is trivial and then use
(\ref{gi}) and (\ref{sym}). We take $[C\cup_q E]$ to be an
arbitrary curve from $R$, where $E$ is an elliptic curve. Using
that limit $\mathfrak g^{3i+4}_{9i+12}$ on $C\cup_q E$ are in
$1:1$ correspondence with linear series $L\in W^{3i+4}_{9i+12}(C)$
having a cusp at $q$ (this being a statement that holds
independent of $E$) and that $\G_{0, b | \sigma^*(\Delta_1^0)}$
consists on each fibre of sections of the genus $g-1$ aspect of
the limit $\mathfrak g^{3i+4}_{9i+12}$, the claim now follows.

Now we  determine $A, B_0$ and $B_1$ explicitly.  We fix a general
pointed curve $(C, q)\in \cM_{6i+9, 1}$ and construct the test
curves $C^1\subset \Delta_1$ and $C^0\subset \Delta_0$. Using the
notation from Proposition \ref{limitlin}, we get that
$\sigma^*(C^0)\cdot  c_1(\G_{i, 2}-\H_{i, 2}) =c_1(\G_{i, 2
|Y})-c_1(\H_{i, 2 |Y})$   and
 $\sigma^*(C^1) \cdot c_1(\G_{i, 2}-\H_{i, 2})=c_1(\G_{i, 2 |X})-c_1(\H_{i, 2 |X})$
 (the other component $T$ of $\sigma^*(C^1)$ does not appear
 because $\G_{0, b |T}$ is trivial for all $b\geq 1$).
On the other hand $$C^0 \cdot \sigma_*(c_1(\G_{i, 2}-\H_{i, 2}))=(12i+18)B_0-B_1 \mbox{ and }
C^1\cdot \sigma_*(c_1(\G_{i, 2}-\H_{i, 2}))=(12i+16)B_1,$$ while we already know that $A-12B_0+B_1=0.$ Using
Propositions \ref{x}, \ref{y}, \ref{has}, the expressions for
$[X]$ and $[Y]$ as well as the well-known formula $\theta^j
x^{k-2-j}=g!/(g-j)!$ on $C_{k-2}$, we get a linear system of $3$
equations in $A, B_0$ and $B_1$ which leads to the stated formulas
for the first three coefficients.

To compute $B_j$ for $2\leq j\leq 3i+5$, we fix general curves
$[C]\in \cM_{6i+10-j}$ and $[D, y]\in \cM_{j, 1}$ and consider the
test curve $C^j\subset \Delta_j$. Then on one hand we have that $2(6i+9-j)B_j=\sigma^*(C^j)\cdot c_1(\G_{i, 2}-\H_{i, 2})$, on the other hand
$$ \sigma^*(C^j)\cdot
c_1(\G_{i,2}-\H_{i, 2})=c_1(\G_{i, 2}-\H_{i, 2})\cdot
\sum_{\alpha=0}^{[j/2]} \Bigl((a(j, \alpha) X_{j, \alpha}+
b(6i+10-j, 3i+5+\alpha-j)Y_{j, \alpha}^{''}\Bigr).$$ Using
(\ref{diferenta}) together with Proposition \ref{secondy}, we get
that
$$c_1(\G_{i, 2}-\H_{i, 2})\cdot \sum_{\alpha=0}^{[j/2]}b(6i+10-j, 3i+5+\alpha-j) Y_{j, \alpha}^{''}
=$$
$$=\frac{(8i^2+19i+12)j! (6i+10-j)! (3i+2)!}{i! (2i+4)!}\Bigl(\sum_{\alpha=0}^j
\frac{(j-2\alpha-1)(j-2\alpha)^2(j-2\alpha+1)}{\alpha!\ (j-\alpha)!\
(3i+5-\alpha)!\ (3i+5+\alpha-j)!}\Bigr),$$ which can be computed by
Maple, whereas $c_1(\G_{i, 2}-\H_{i, 2})\cdot
(\sum_{\alpha=0}^{[j/2]} a(j, \alpha) X_{j, \alpha})$ can be
computed using (\ref{diferenta}) and Propositions \ref{xja} and
\ref{xjb} (note in particular that $\aleph$ has to be multiplied by the coeffcient of $c_1(\G_{0, 1})$ in (\ref{diferenta})). The rest now follows and checking that $b_j\geq b_0$ for
$j\geq 1$ is elementary.
\end{proof}

\begin{remark}
If we retain the notation of Theorem \ref{formula}, we see that if we fix $j\geq 1$ and let $i$ vary, we have that $\mbox{lim}_{i\rightarrow \infty} \frac{b_j}{b_0}=6j$. This is the same asymptotical estimate as in the case of the classical Brill-Noether divisors where $\frac{b_j}{b_0}=\frac{6j(g-j)}{g+3}$ for all $j\geq 1$ (cf. \cite{EH3}).
\end{remark}

\begin{theorem}\label{as}
Assume that vector bundle morphism $\phi:\H_{i, 2}\rightarrow
\G_{i, 2}$ is non-degenerate at a general point from
$\sigma^{-1}(\Delta_1^0)$. Then $\zz_{6i+10, i}$ is a divisor on
$\mm_{6i+10}$ and
$$6<s(\zz_{6i+10,i})=\frac{3(4i+7)(6i^2+19i+12)}{(12i^2+31i+18)(i+2)}<6+\frac{12}{g+1}.$$
\end{theorem}

\begin{remark} Since $\sigma^{-1}(\Delta_1^0)$ is irreducible, to
verify the assumption made in Theorem \ref{as} it suffices to show
that if $C\subset \PP^{3i+4}$ is a general $1$-cuspidal curve with
 $p_a(C)=6i+10$ and $\mbox{deg}(C)=9i+12$, then $C$ satisfies property$(N_i)$.
We checked this for $g=16, 22$ (cf. Remarks
\ref{cusp16} and \ref{cusp22}), while the case $g=10$ was treated
in \cite{FP}. The assumption in Theorem  \ref{as} implies that
$\phi$ is non-degenerate over a general point from
$\sigma^{-1}(\Delta_0^0)$ and also over a general point from
$\sigma^{-1}(\cM_g^0)$, which is equivalent to $\mathcal{Z}_{g,
i}$ being a divisor on $\cM_g^0$. Therefore the assumption is
slightly stronger than Conjecture \ref{stratum}.
\end{remark}

\noindent \emph{Proof of Theorem \ref{as}.} We denote by $Z(\phi)$
the degeneracy locus of $\phi$, thus $\sigma_*(Z(\phi))$ is a
divisor on $\mm_{6i+10, i}$ such that $\sigma_*(Z(\phi))=\zz_{g,
i}+\sum_{j=0}^{3i+5} d_j \Delta_j$, for certain coefficients
$d_j\geq 0$. Using the Remark above we obtain that $d_0=d_1=0$,
hence $s(\zz_{6i+10, i})=s(\sigma_*(c_1(\G_{i, 2}-\H_{i, 2})))$ and
the rest follows from Theorem \ref{formula}.
\hfill $\Box$

\begin{remark} In the simplest case $i=0$, we can compare the formula for
$[\zz_{6i+10, i}]$ with our findings in \cite{FP}. Theorem \ref{formula} gives
the formula
$$\sigma_*(c_1(\G_{0, 2}-\H_{0, 2}))\equiv 42(7 \lambda- \delta_0- 5 \delta_1-\frac{4}{3}\delta_2-7\delta_3-
\frac{44}{7}\delta_4-\frac{685}{42} \delta_5),$$ whereas Theorem 1.6
from \cite{FP} says that if $\mathcal{K}_{10}$ is the divisor on
$\cM_{10}$ of curves lying on a $K3$ surface then
$$\kk_{10}=7\lambda-\delta_0-5\delta_1-9\delta_2-12\delta_3-14\delta_4-b_5\delta_5,$$
where $b_5\geq 6$. Since we have also established the set-theoretic
equality $\mathcal{Z}_{10, 0}=\mathcal{K}_{10}$ (cf. \cite{FP},
Theorem 1.7), it follows that there is a scheme-theoretic equality
$\mathcal{Z}_{10, 0}=42 \mathcal{K}_{10}$. Here $42$ is the number
of pencils $\mathfrak g^1_6$ on a general curve of genus $10$, and
its appearance in the formula of $[\zz_{10, 0}]$ has a clear
geometric meaning: if a Brill-Noether general curve $[C]\in
\cM_{10}$  fails $(N_0)$ for {\emph{one}} linear system $\mathfrak
g^4_{12}=K_C(-\mathfrak g^1_6)$, then it fails $(N_0)$ for all $42$
linear systems $\mathfrak g^4_{12}$ it possesses, that is, the map
$\sigma:\mathcal{U}_{10, 0}\rightarrow \K_{10, 0}$ is $42:1$.
Moreover, the vector bundle morphism $\H_{0, 2}\rightarrow \G_{0,
2}$ is degenerate along each of the boundary divisors $\Delta_2,
\ldots, \Delta_5$. This situation presumably extends to higher
genera as well hence the main thrust of Theorem \ref{formula} is
that it computes $s(\zz_{6i+10, i})$.

\end{remark}

To finish the proof of Theorems \ref{g16} and \ref{g22} we
specialize to the cases $g=16$ and $22$ and note that the assumption
made in Theorem  \ref{as} is satisfied in these situations (cf.
Theorems \ref{gen16} and \ref{gen22}).

\section{An effective divisor on $\mm_{14, 1}$}

In this section we describe how to construct effective divisors on
$\mm_{g, n}$ using syzygy type conditions for pointed curves. We
treat only one example. We denote by $\zz_{14, 0}^1$ the closure
in $\mm_{14, 1}$ of the locus $\mathcal{Z}_{14, 0}^1$ of smooth
pointed curves $[C, p]\in \cM_{14, 1}$ for which there exists a
linear series $L\in W^6_{18}(C)$ such that the map
$$\mu_{L, p}:\mbox{Sym}^2 H^0(C, L(-p))\rightarrow H^0(C, L^{\otimes
2}(-2p))$$ is not an isomorphism. Just like the in the case of the
loci $\mathcal{Z}_{g, i}$ on $\cM_g$, the locus $\mathcal{Z}_{14,
0}^1$ can be naturally viewed as the pushforward of the degeneracy
locus of a morphism between two vector bundles of the same rank
$21$ over $\mathfrak G^6_{18}\times_{\cM_{14}} \cM_{14, 1}$.

We are going to show that $\zz_{14, 0}^1$ is a divisor and in
order to compute its class we need some preparations. Recall that
for $g\geq 3$, the group $\mbox{Pic}(\mm_{g, 1})$ is freely
generated by $\lambda$, the tautological class $\psi$ and the
boundary classes $\delta_i=[\Delta_i]$ with $0\leq i\leq g-1$,
where for $i\geq 1$, the generic point of  $\Delta_i$ is a union
of two smooth curves of genus $i$ and $g-i$ meeting at a point,
the marked point lying on the genus $i$ component. We denote by
$\pi:\mm_{g, 1}\rightarrow \mm_g$ the natural forgetful map. We
can then write the class of $\zz_{14, 0}^1$ on $\mm_{14, 1}$ as
$$\zz_{14, 0}^1\equiv a \ \lambda+c \ \psi-b_0\ \delta_0-b_1 \ \delta_1-\cdots
-b_{13} \ \delta_{13}.$$

As before, we will determine the relevant coefficients in the
expression of $[\zz_{14, 0}^1]$ by intersecting our divisor with
various test curves.

\begin{proposition} Let $C$ be a general curve of genus $14$ and
$\widetilde{C}=\pi^{-1}([C])\subset \mm_{14, 1}$, the test curve
obtained by letting the marked point vary along $C$. Then
$\widetilde{C}\cdot \zz_{14, 0}^1=12012$, hence
$c=\widetilde{C}\cdot \zz_{14, 0}^1/(2g-2)=462$.
\end{proposition}
\begin{proof} Let us fix a linear series $L\in W^6_{18}(C)$. We
count the number of points $p\in C$, for which the multiplication
map $\mu_{L, p}$ not injective. If $p_1:C\times C\rightarrow C$
and $p_2:C\times C\rightarrow C$ are the two projections, we
define the vector bundles
$$\E:=(p_{2})_*(p_{1}^*(L)\otimes \OO_{C\times C}(-\Delta))\mbox{
and } \F:=(p_{2})_*(p_1^*(L^{\otimes 2})\otimes \OO_{C\times
C}(-2\Delta)).$$ There is a natural multiplication map
$\mu_L:\mbox{Sym}^2(\E)\rightarrow \F$, and the cardinality of its
degeneracy locus is just $c_1(\F)-c_1(\mbox{Sym}^2(\E))$. A simple
calculation shows that $c_1(\E)=-4x-\theta=-18$, hence
$c_1(\mbox{Sym}^2(\E))=-126$, while $c_1(\F)=-98$. We obtain that
$\mu_{L, p}$ is not an isomorphism for precisely $28$ points $p\in
C$. Since $C$ has $429$ linear series $\mathfrak g^6_{18}$ (cf.
\cite{ACGH}), we obtain $c=(\tilde{C}\cdot \zz_{14,
0}^1)/26=429\cdot 28/26=462$.
\end{proof}

\noindent For more relations among the coefficients of $[\zz_{14,
0}^1]$ we define the map $\nu:\mm_{1, 2}\rightarrow \mm_{g, 1}$
obtained by attaching to each $2$-pointed elliptic curve $[E, q,
p]$ a fixed general $1$-pointed curve $[C, q]\in \cM_{g-1}$ (the
point of attachment being $q$). One has the pullback formulas
$$\nu^*(\lambda)=\lambda, \ \nu^*(\psi)=\psi_p, \mbox{ }\nu^*(\delta_0)=\delta_0, \mbox{
}\nu^*(\delta_1)=-\psi_q\mbox{ and
}\nu^*(\delta_{g-1})=\delta_{qp},$$ where $\psi_p$ and $\psi_q$
are the tautological classes corresponding to the marked points
$p$ and $q$, while $\delta_{qp}$ is the boundary component of
curves having a rational tail containing both $q$ and $p$. On
$\mm_{1, 2}$ these classes are not independent and we have the
relations $\psi_q=\psi_p$, $\lambda=\psi_p-\delta_{qp}$ and
$\delta_0=12(\psi_p-\delta_{qp})$ (see e.g. \cite{AC2},
Proposition 1.9).

\begin{proposition}\label{m12}
If $\nu:\mm_{1, 2}\rightarrow \mm_{14, 1}$ is as above, then
$\nu^*(\zz_{14, 0}^1)=\emptyset$. It follows that we have the
relations
$$a-12b_0+b_{13}=0\mbox{ and }c+b_1=b_{13}.$$
\end{proposition}
\begin{proof}
We assume that  $[X=C\cup_q E, p\in E]\in \zz_{14, 0}^1$. Then
there exists a limit linear series $\mathfrak g^6_{18}$ on $X$
determined by its aspects $L_C\in W^6_{18}(C)$ and $L_E\in
W^6_{18}(E)$, together with compatible elements
$$\rho_C\in \mbox{Ker}\{\mu_{L_C}:\mbox{Sym}^2 H^0(L_C)\rightarrow
H^0(L_C^{\otimes 2})\}, \mbox{ }\rho_E\in
\mbox{Ker}\{\mu_{L_E}:\mbox{Sym}^2 H^0(L_E)\rightarrow
H^0(L_E^{\otimes 2})\},$$ satisfying the inequality
$$\mbox{ord}_q(\rho_C)+\mbox{ord}_q(\rho_E)\geq
\mbox{deg}(L_C)+\mbox{deg}(L_E)=36$$ and such that $\rho_E\in
\mbox{Sym}^2 H^0(L_E(-p))$ (see \cite{FP}, Section 4, for how to
study multiplication maps in the context of limit linear series).
Because $(C, q)\in \mm_{13, 1}$ is general, the vanishing sequence
of $L_C$ at $q$ equals $(0, 2, 3, 4, 5, 6, 7)$, the vanishing
sequence of $L_E$ at $q$ is $(11, 12, 13, 14, 15, 16, 18)$ and
finally, the vanishing sequence of $L_E$ at $p$ is either $(0, 1,
2, 3, 4, 5, 6)$ or $(0, 1, 2, 3, 4, 5, 7)$, depending on whether
$p-q\in \mbox{Pic}^0(E)$ is  a $7$-torsion class or not.

We claim that $\mbox{ord}_q(\rho_E)\leq 29(=13+16=14+15)$. Indeed,
otherwise $\mbox{ord}_q(\rho_E)\geq 30(=14+16=15+15)$, and since
$\rho_E\in \mbox{Sym}^2 H^0(L_E(-p))$, after subtracting base
points  $\rho_E$ becomes a $\neq 0$ element in the kernel of the
map $\mbox{Sym}^2 H^0(N)\rightarrow H^0(N^{\otimes 2})$, where
$N=L_E(-p-14q)\in \mbox{Pic}^3(E)$. This is obviously impossible.
Therefore $\mbox{ord}_q(\rho_E)\leq 29$, so by compatibility,
$\mbox{ord}_q(\rho_C)\geq 7(=2+5=3+4)$. We now show that when
$(C,p)\in \cM_{13, 1}$ is chosen generically, there can be no such
element $\rho_C$ which leads to a contradiction.

\noindent \textbf{Claim:} Suppose $\sigma_0, \sigma_2, \sigma_3,
\sigma_4, \sigma_5, \sigma_6, \sigma_7$ is a basis of $H^0(L_C)$
adapted to the point $q$ in the sense that
$\mbox{ord}_q(\sigma_i)=i$. If $W(q, L_C)\subset \mbox{Sym}^2
H^0(L_C)$ denotes the $17$-dimensional subspace spanned by the
elements $\sigma_i\cdot \sigma_j$ with $4\leq i\leq j\leq 7$,\
$\sigma_3\cdot \sigma_j$ for $j\geq 4$ and $\sigma_2\cdot
\sigma_j$ for $j\geq 5$, then the restriction of the
multiplication map $W(q, L_C)\rightarrow H^0(L_C^2(-7q))$ is an
isomorphism. (Note that $W(q, L_C)$ does not depend on the chosen
basis $\{\sigma_i\}$).

The proof of this claim is similar to the proof of Theorem 5.1 in
\cite{FP}. It is enough to construct a single $1$-cuspidal curve
$X\subset \PP^6$ with $p_a(X)=14$ and $\deg(X)=18$, such that if
$\nu:C\rightarrow X$ is the normalization of $X$ and $q\in C$ is
the inverse image of the cusp, then $X$  does not lie on any
quadric contained in $W\bigl(q, \nu^*(\OO_X(1))\bigr)$. We
construct the following cuspidal curve: define $\Gamma \subset
\PP^6$ to be the image of the map $t\stackrel{f}\mapsto [1, t^2,
t^3, t^4, t^5, t^6, t^7]$, then choose a general hyperplane
$H\subset \PP^6$ which intersects $\Gamma$ in distinct points
$p_1, \ldots, p_7$. Take $D\subset H$ to be a general smooth curve
of genus $g(D)=7$ and $\mbox{deg}(D)=11$ which passes through
$p_1, \ldots ,p_7$. Then $X:=\Gamma\cup D$ is a  curve of
arithmetic genus $14$ and degree $21$ having a cusp at the point
$q=f(0)\in \Gamma$. The quadrics in $W\bigl(q,
\nu^*(\OO_X(1))\bigr)$ can of course be written down explicitly
and to show that $D$ can be chosen such that it is not contained
in any quadric from $W\bigl(q, \nu^*(\OO_X(1))\bigr)$ amounts to a
simple counting argument.
\end{proof}

\begin{remark} The claim we have just proved also shows that $\zz_{14,
0}^1$ is a divisor on $\mm_{14, 1}$. Alternatively, this can be
proved in the same way as Theorem \ref{gen16}.
\end{remark}

We now determine the coefficient $b_1$:
\begin{proposition} Consider general curves $[E, p, q]\in \cM_{1,
2}$ and $[C]\in \cM_{13}$, and denote by $C^2\subset \mm_{14, 1}$
the test curve consisting of points $\{[C\cup_q E, p\in E]\}_{\{q
\in C\}}$. Then $C^2\cdot \zz_{14, 0}^1=133848$, therefore
$b_1=C^2\cdot \zz_{14, 0}^1/24=5577.$
\end{proposition}
\begin{proof} Throughout the proof we will use the notations introduced
in Proposition \ref{x}. From the proof of Proposition \ref{m12} it
is clear that the intersection number $C^2\cdot \zz_{14, 0}^1$
equals the number of pairs $(q, D)\in X\subset C\times C_6$ such
that the map
$$W(q, K_C(-D))\rightarrow H^0\bigr(K_C^{\otimes
2}(-2D-7q)\bigl)$$ is not an isomorphism. To compute this number
we construct the rank $17$ vector bundle $W$ over $X$ having fibre
$W(q, D)=W(q, K_C(-D))$ over each point $(q, D)\in X$. The curve
$C$ being general, the vanishing sequence $a^{K_C(-D)}(q)$ will be
generically equal to $(0, 2, 3, 4, 5, 6, 7)$, while at a finite
number of points $(q, D)\in X$ we will have that
$a^{K_C(-D)}(q)=(0, 2, 3, 4, 5, 6, 8)$. In order to compute
$c_1(W)$ we note that $W$ has a subbundle $W_1\subset W$ which
fits into two exact sequences:
$$0\longrightarrow \mbox{Sym}^2 u_*\bigl(v^*(\cM)\otimes
\OO(-4\Delta)\bigr)_{|X}\longrightarrow W_1 \rightarrow
u_*\bigl(v^*(\cM)\otimes \OO(-4\Delta)\bigl)_{|X}\otimes
\mathcal{P}\longrightarrow 0,$$ where $\mathcal{P}=u_*\bigl(
v^*(\cM)\otimes \I_{\Delta}^3/{\I_{\Delta}^4}\bigr)_{|X}$, and
$$0\longrightarrow W_1\longrightarrow W\longrightarrow
u_*\bigl(v^*(\cM)\otimes \OO(-5\Delta)\bigl)_{| X}\otimes
u_*\bigl(v^*(\cM)\otimes
\I_{\Delta}^2/{\I_{\Delta}^3}\bigr)_{|X}\longrightarrow 0.$$ Next,
we consider the multiplication map $W\rightarrow
u_*\bigl(v^*(\cM^{\otimes 2})\otimes \OO(-7\Delta)\bigr)_{| X}$
whose degeneration locus we want to compute. The intersection
number $C^2\cdot \zz_{14, 0}^1$ is  equal to
$$\Bigl(c_1\bigl(u_*(v^*(\cM^{\otimes 2})\otimes
\OO(-7\Delta))\bigr)-c_1(W)\Bigr)\cdot
[X]=(3\theta-7x-\gamma+18\eta)\cdot [X]=133848,$$ as it turns out
after a short calculation. Here we have used that $c_1(W)$ can be
computed from the two exact sequences involving $W$ and $W_1$,
while $$c_1\bigl(u_*(v^*(\cM^{\otimes 2})\otimes
\OO(-7\Delta))\bigr)=c_1(\mathcal{F}_2)-c_1(J_6(\cM^{\otimes
2}))=-4\theta-34x+14\gamma-756\eta.$$
\end{proof}

Since $b_1=5577$ we now obtain that $b_{13}=b_1+c=6039$. To
compute the coefficient $b_0$ (and thus the $\lambda$-coefficient
$a$) we use our last test curve:

\begin{proposition}
Let $(C, q, p)\in \cM_{13, 2}$ be a general $2$-pointed curve. We
denote by $C^3\subset \mm_{14, 1}$ the family  consisting of
curves $\{[C/y\sim q, p]\}_{\{y\in C\}}$. Then $C^3\cdot \zz_{14,
0}^1=c+26b_0-b_{13}=24453$. It follows that $b_0=1155$ and
$a=7821$.
\end{proposition}

\begin{proof} We retain the notations introduced in Proposition
\ref{y}. We construct a vector bundle map
$\mbox{Sym}^2(\E)\rightarrow \F$ over the curve $Y\subset C\times
C_6$, where $\E=u_*\bigl(v^*(\cM)\otimes \OO(-\Gamma_p)\bigr)$ is
the bundle with fibre $\E(y, D)=H^0\bigl(K_C(-D-p)\bigr)$, while
$\F$ is the bundle with fibre
$$\F(y, D)=H^0\bigl(K_C^{\otimes 2}(-2D-2p-y-q)\bigr)\oplus \mathbb C
\cdot t^2\subset H^0\bigl(K_C^{\otimes ^2}(-2D-2p)\bigr),$$ where
$H^0\bigl(K_C(-D-p)\bigr)/H^0\bigl(K_C(-D-p-y-q)\bigr)=\mathbb C
\cdot t$, for every $(y, D)\in Y$.

It is easy to show that $c_1(\E)=-5x-\theta$, hence
$c_1(\mbox{Sym}^2(\E))=-35x-7\theta$. The class  $c_1(\F)$ can be
computed from the exact sequence
$$0\longrightarrow u_*\bigl(v^*(\cM^{\otimes 2})\otimes
\OO(-2\Gamma_p-\Gamma_q-\Delta)\bigr)_{| Y}\longrightarrow
\F\longrightarrow \bigl(\cM^{\otimes 2}\otimes
\OO(-2\Gamma_p)\bigr)_{| Y}(-A)\longrightarrow 0,$$ where $A=Y\cap
\pi_1^{-1}(p)$ is the effective divisor on the curve $Y$
consisting of all points $(p, D)$ such that $h^0(p+q+D)\geq 2$.
Using the formula for $[Y]$ (cf. Proposition \ref{y}), we get that
$|A|=429$. We also compute that $c_1\bigl(u_*(v^*(\cM^{\otimes
2})\otimes
\OO(-2\Gamma_p-\Gamma_q-\Delta))\bigr)=-4\theta-40x-33\eta+2\gamma$
and of course $c_1(\cM^{\otimes 2}\otimes \OO(-2\Gamma_p))=2(17
\eta-\gamma-x)$. Therefore we can write that
$$C^3\cdot \zz_{14, 0}^1=\bigl(c_1(\F)-c_1(\mbox{Sym}^2(\E))\bigr)\cdot
[Y]-|A|=(-7x+3\theta+\eta)\cdot [Y]-|A|=24453.$$
\end{proof}

We have thus far determined the coefficients $a, c, b_0, b_1,
b_{13}$ in the expansion of $[\zz_{14, 0}^1]$. This is already
enough to conclude that $[\zz_{14, 0}^1]$ lies outside the cone of
$\mbox{Pic}(\mm_{14, 1})$ spanned by pullbacks of effective
divisors from $\mm_{14}$, Brill-Noether divisors on $\mm_{14, 1}$
and boundary divisors (see the discussion after Theorem \ref{m14}
for the relevance of this result). To get a bound on the remaining
coefficients $b_j$ for $2\leq j\leq 12$, we use a variant of
Theorem 1.1 from \cite{FP}. The boundary divisor $\Delta_j\subset
\mm_{14, 1}$ with $3\leq j\leq 12, j\neq 4$ is filled-up by
pencils $R_j$ obtained by attaching to a fixed $2$-pointed curve
$[B, p, q]\in \cM_{j, 2}$ a variable $1$-pointed curve $[C, q]\in
\mm_{14-j, 1}$ moving in a Lefschetz pencils of curves of genus
$14-j$ sitting on a fixed $K3$ surface. Deformations of $R_j$
cover $\Delta_j$ for $j\neq 2, 4$, hence we have that $R_j\cdot
\zz_{14, 0}^1\geq 0$. Since one also has the relations (see
\cite{FP}, Lemma 2.4)
$$R_j\cdot \lambda=15-j, R_j\cdot \delta_j=-1, R_j\cdot
\delta_0=6(17-j), R_j\cdot \psi=0 \mbox{ and }R_j\cdot
\delta_i=0\mbox{ for } i\neq 0, j,$$ we immediately get the
estimate $b_j\geq 15+27j$ for all $3\leq j\leq 12$, $j\neq 4$. To
obtain the bounds $b_2\geq 325$ and $b_4\geq 271$ we use similar
pencils filling up $\Delta_2$ and $\Delta_4$ respectively. We skip
these details. This completes the proof of Theorem \ref{m14}.

\section{The Kodaira dimension of $\mm_{g, n}$}
In this last section we use the effective divisors $\zz_{16,1 }$
and $\zz_{22, 2}$ to improve Logan's results about which moduli
spaces $\mm_{g, n}$ are of general type. For a general reference
about $\mbox{Pic}(\mm_{g, n})$ we refer to \cite{AC2} and
\cite{Lo}. For each $1\leq i\leq n$ we denote by $\pi_i:\mm_{g,
n}\rightarrow \mm_{g, 1}$ the morphism forgetting all marked
points except the one labelled by $i$ and we also consider the map
$\pi:\mm_{g, n}\rightarrow \mm_g$ which forgets all marked points.
We recall that the canonical class of $\mm_{g,n}$ is given by the
formula
$$K_{\mm_{g, n}}=13\lambda-2\delta_0+\sum_{i=1}^n
\psi_i-2\sum_{i\geq 0, \ S}\delta_{i:S}-\sum_{S} \delta_{1: S}.$$
Here $\psi_i$ is the tautological class corresponding to the
$i$-th marked point, while $\delta_{i: S}$ with $i\geq 0$ and
$\emptyset \neq S\subset \{1, \ldots, n\}$ denotes the class of
the boundary divisor with generic point being a union of two
curves of genus $i$ and $g-i$ such that the marked points on the
genus $i$ component are precisely those labelled by $S$.

\noindent \emph{Proof of Theorem \ref{mgn}.} We start with the
case $g=22$ and we show that $\mm_{22, 2}$ is of general type. On
$\mm_{22, 2}$ we consider the averaged pullback of the Weierstrass
divisor
$$W_{12}:=\pi_1^*(\overline{\mathcal{W}})+\pi_2^*(\overline{\mathcal{W}})\equiv
-2\ \lambda+\frac{g(g+1)}{2}(\psi_1+\psi_2)-\sum_{i\geq 0, S\neq
\emptyset} c_{i: S}\ \delta_{i:S},$$ where $c_{i: S}\geq 0$, and
significantly, the coefficient of $\delta_0$ is $0$.  By Theorem
\ref{g22}, we have another effective  class, namely
$[\pi^*(\zz_{22, 2})]=c(\frac{1665}{256}
\lambda-\delta_0-\cdots)$, where $c>0$. One can easily check that
$K_{\mm_{22, 2}}$ can be written as a positive combination of
$[W_{12}], [\pi^*(\zz_{22, 2})], \psi_1+\psi_2$ and some other
boundary classes. Since $\psi_1+\psi_2$ is big and nef, it follows
that $\mm_{22, 2}$ is of general type.

When $g=21$  we use the maps $\chi_{i, j}:\mm_{21, 5}\rightarrow
\mm_{22}$ for $1\leq i<j\leq 5$, where $\chi_{i, j}$ identifies
the marked points labelled by $i$ and $j$ and forgets those
labelled by $\{i, j\}^c$. The $\QQ$-divisor class
$$\sum_{i<j} \chi_{i, j}^*(\zz_{22, 2})\equiv c(\frac{1665}{256}\lambda-\delta_0+\frac{2}{5} \sum_{i=1}^5
\psi_i-\cdots), \ \mbox{ where }c>0, $$ is obviously effective on
$\mm_{21, 5}$. Since the $\QQ$-class $-\lambda+11\sum_{i=1}^5
\psi_i-0\cdot \delta_0-\cdots$, is also effective (cf. \cite{Lo},
Theorem 5.4 - we have retained only the coefficients that play a
role in our argument), once again we see that $K_{\mm_{21, 5}}$
can be written as the sum of an effective divisor and a positive
multiple of $\sum_{i=1}^5 \psi_i$.

Finally we settle the case $g=16$: we adapt Theorem 5.4 from
\cite{Lo} to conclude that the $\QQ$-class
$-\lambda+\frac{23}{9}(\sum_{i=1}^9 \psi_i)-0\cdot
\delta_0-\cdots$, is effective on $\mm_{16, 9}$ (precisely, this
is the class of the $S_9$-orbit of the closure in $\mm_{16, 9}$ of
the effective divisor $D$ on $\cM_{16, 9}$ consisting of points
$[C, p_1,\ldots, p_9]$ such that $h^0(C,
2p_1+\cdots+2p_7+p_8+p_9)\geq 2$ and we have explicitly indicated
that the coefficient of $\delta_0$ is $0$ and retained only the
coefficients that are significant for this calculation).  The
class $[\pi^*(\zz_{16,1})]=c(\frac{407}{61} \
\lambda-\delta_0-\cdots)$, where $c>0$, is also effective and one
writes $K_{\mm_{16, 9}}$ as a positive sum of these two effective
classes, boundary classes and the big and nef class $\sum_{i=1}^9
\psi_i$.\hfill $\Box$

\begin{remark} Theorem 6.3 from \cite{Lo} claims that $\mm_{22,
4}$ is of general type but the numerical argument used in the
proof seems to be incorrect.
\end{remark}

\end{document}